\documentclass[12pt,draft]{amsart}

\usepackage[cp1251]{inputenc}
\usepackage[T2A]{fontenc}
\usepackage[english]{babel}
\usepackage{amsmath,amsfonts,amssymb}
\usepackage{geometry}
\newtheorem{theorem}{Theorem}[section]

\newtheorem{lemma}{Lemma}[section]

\theoremstyle{remark}
\newtheorem{remark}{Remark}[section]

\theoremstyle{definition}
\newtheorem{definition}{Definition}[section]

\author[V. A. Mikhailets, A. A. Murach]{Vladimir A. Mikhailets, Alexandr A. Murach}

\title[Interpolation and refined scale]{Interpolation with a
function parameter and \\ refined scale of spaces}

\address{Institute of Mathematics NAS of Ukraine,
Tereshchenkivska str. 3, Kyiv, Ukraine, 01601}

\email{mikhailets@imath.kiev.ua}

\address{Institute of Mathematics NAS of Ukraine,
Tereshchenkivska str. 3, Kyiv, Ukraine, 01601; \indent Chernigiv
State Technological University, Shevchenka str. 95, Chernigiv
14027, Ukraine}

\email{murach@imath.kiev.ua}


\date{07/12/2007}


\begin{document}

\maketitle

\begin{abstract}
The interpolation of couples of separable Hilbert spaces with a
function parameter is studied. The main properties of the classic
interpolation are proved. Some applications to the interpolation
of isotropic H\"ormander spaces over a closed manifold are given.
\end{abstract}

\section{Introduction}
In this paper we study the interpolation of couples of separable
Hilbert spaces with a functional parameter. We generalize the
classical theorems on interpolation with a power parameter of
order $\theta\in(0,1)$ to the maximal class of functions.

As an application, we consider the interpolation of isotropic
H\"ormander spaces over a closed manifold
$$
H^{s,\varphi}:=H_{2}^{\langle\cdot\rangle^{s}\,\varphi(\langle\cdot\rangle)},
\quad\langle\xi\rangle:=\bigl(1+|\xi|^{2}\bigr)^{1/2}. \leqno(1.1)
$$
Here, $s\in\mathbb{R}$ and $\varphi$ is a functional parameter
slowly varying at $+\infty$ in Karamata's sense. In particular,
every standard function
$$
\varphi(t)=(\log t)^{r_{1}}(\log\log
t)^{r_{2}}\ldots(\log\ldots\log
t)^{r_{n}},\quad\{r_{1},r_{2},\ldots,r_{n}\}\subset\mathbb{R},\;n\in\mathbb{N},
$$
is admissible. This scale was introduced and investigated by the
authors in \cite{MM2a, MM2b}. It contains the Sobolev scale
$\{H^{s}\}\equiv\{H^{s,1}\}$ and is attached to it by the number
parameter $s$ and being considerably finer.

Spaces of form (1.1) arise naturally in different spectral
problems: convergence of spectral expansions of self-adjoint
elliptic operators almost everywhere, in the norm of the spaces
$L_{p}$ with $p>2$ or $C$ (see survey \cite{AIN76}); spectral
asymptotics of general self-adjoint elliptic operators in a
bounded domain, the Weyl formula, a sharp estimate of the
remainder in it (see \cite{Mikh82, Mikh89}) and others. They may
be expected to be useful in other "fine"\, questions.  Due to
their interpolation properties, the spaces $H^{s,\varphi}$ occupy
a special position among the spaces of a generalized smoothness
which are actively investigated and used today (see survey
\cite{KalLiz87}, recent articles \cite{HarMou04, FarLeo06} and the
bibliography given there).

One of the main results of the article is a description of the
refined scale by means of regularly varying functions of a
positive elliptic pseudodifferential operator.

The related questions were studied in \cite{Merucci84, CobFern86}
and by the authors in [11--20].

\section{An interpolation with a function parameter}

\subsection{A definition of the interpolation.}

\begin{definition}
An ordered couple $[X_{0},X_{1}]$ of complex Hilbert spaces
$X_{0}$ and $X_{1}$ is called \textit{admissible} if the spaces
$X_{0}$, $X_{1}$ are separable and the continuous imbedding
$X_{1}\hookrightarrow X_{0}$ holds.
\end{definition}

\vspace{0.3 cm}

Let an admissible couple $X=[X_{\,0},X_{1}]$ of Hilbert spaces be
given. It is known \cite[Ch. 1, Sec. 2.1]{LM71}, \cite[Ch. IV,
Sec. 9.1]{Kr72} that for this couple $X$ there exists the
isometric isomorphism $J:X_{1}\leftrightarrow X_{\,0}$ such that
$J$ is a self-adjoint positive operator on the space $X_{\,0}$
with the domain $X_{1}$. The operator $J$ is called a
\textit{generating} one for the couple $X$. This operator is
uniquely determined by the couple\;$X$. Indeed, assume that
$J_{1}$ is also a generating operator for the couple $X$. Then the
operators $J$ and $J_{1}$ are metrically equal, that is
$\|Ju\|_{X_{\,0}}=\|u\|_{X_{1}}=\|J_{1}u\|_{X_{\,0}}$ for any
$u\in X_{1}$. Moreover, these operators are positive. Hence, they
are equal.

We denote by $\mathcal{B}$ the set of all functions
$\psi:(0,+\infty)\rightarrow(0,+\infty)$ such that
\begin{description}
  \item[a)] $\psi$ is Borel measurable on the semiaxis $(0,+\infty)$;
  \item[b)] $\psi$ is bounded on each closed interval $[a,b]$, where
$0<a<b<+\infty$;
  \item[c)] $1/\psi$ is bounded on each set $[r,+\infty)$, where
  $r>0$.
\end{description}

Let $\psi\in\mathcal{B}$. Generally, the unbounded operator
$\psi(J)$ is defined in the space $X_{0}$ as a function of $J$. We
denote by $[X_{0},X_{1}]_{\psi}$ or, simply, by $X_{\psi}$ the
domain of the operator $\psi(J)$, equipped with the inner product
$(u,v)_{X_{\psi}}:=(\psi(J)u,\psi(J)v)_{X_{0}}$ and the
corresponding norm $\|\,u\,\|_{X_{\psi}}=(u,u)_{X_{\psi}}^{1/2}$.

The space $X_{\psi}$ is the Hilbert separable one and, moreover,
the continuous dense imbedding $X_{\psi}\hookrightarrow X_{0}$ is
fulfilled. Indeed, we have $\mathrm{Spec}\,J\subseteq[r,+\infty)$
and $\psi(t)\geq c$ for $t\geq r$, where $r,c$ are some positive
numbers. Hence, $\mathrm{Spec}\,\psi(J)\subseteq[c,+\infty)$, that
implies the isometric isomorphism $\psi(J):X_{\psi}\leftrightarrow
X_{0}$. It follows that the space $X_{\psi}$ is complete and
separable as well as that the function $\|\cdot\|_{X_{\psi}}$ is
positive definite, so this is a norm. Next, since the operator
$\psi^{-1}(J)$ is bounded in the space $X_{0}$, a bounded
imbedding operator $I=\psi^{-1}(J)\psi(J):X_{\psi}\rightarrow
X_{0}$ exists. The imbedding $X_{\psi}\hookrightarrow X_{0}$ is
dense because the domain of the operator $\psi(J)$ is a dense
linear manifold in the space $X_{0}$.

Further, it is useful to note the following. Let functions
$\varphi,\psi\in\mathcal{B}$ be such that $\varphi\asymp\psi$ in a
neighborhood of $+\infty$. Then, by the definition of the set
$\mathcal{B}$, we have $\varphi\asymp\psi$ on $\mathrm{Spec}\,J$.
Hence, $X_{\varphi}=X_{\psi}$ up to equivalent norms.

\vspace{0.3 cm}

\begin{definition}
A function $\psi\in\mathcal{B}$
is called an \textit{interpolation parameter} if the following
condition is satisfied for all admissible couples
$X=[X_{0},X_{1}]$, $Y=[Y_{0},Y_{1}]$ of Hilbert spaces and an
arbitrary linear mapping $T$ given on $X_{0}$: \; if the
restriction of the mapping $T$  to the space $X_{j}$ is a bounded
operator $T:X_{j}\rightarrow Y_{j}$ for each $j=0,\,1$, then the
restriction of the mapping $T$ to the space $X_{\psi}$ is also a
bounded operator $T:X_{\psi}\rightarrow Y_{\psi}$.
\end{definition}

Otherwise speaking, $\psi$ is an interpolation parameter if and
only if the mapping $X\mapsto X_{\psi}$ is an interpolation
functor given on the category of all admissible couples $X$ of
Hilbert spaces \cite[Sec. 1.2.2]{Tr80}, \cite[Sec. 2.4]{BL80}. In
the case where $\psi$ is an interpolation parameter, we say that
the space $X_{\psi}$ is obtained by the \textit{interpolation with
the function parameter} $\psi$ of the admissible couple $X$.

Further we investigate the main properties of the mapping
$X\mapsto X_{\psi}$.

\vspace{0.3 cm}

\subsection{Imbeddings of spaces.}

\begin{theorem}
 Let $\psi\in\mathcal{B}$ be an
interpolation parameter and $X=[X_{0},X_{1}]$ be an admissible
couple of Hilbert spaces. Then the continuous dense imbeddings
$X_{1}\hookrightarrow X_{\psi}\hookrightarrow X_{0}$ hold.
\end{theorem}

\begin{proof}
According to Subsection 2.1 it only remains to prove the
continuous dense imbedding $X_{1}\hookrightarrow X_{\psi}$.
Consider two bounded imbedding operators $I:X_{1}\rightarrow
X_{0}$ and $I:X_{1}\rightarrow X_{1}$. Since $\psi$ is an
interpolation parameter, these operators imply the bounded
imbedding operator $I:X_{1}\rightarrow X_{\psi}$. Thus, the
continuous imbedding $X_{1}\hookrightarrow X_{\psi}$ is valid. We
will prove that it is dense. For an arbitrary $u\in X_{\psi}$, we
have $v:=(1+\psi^{2}(J))^{1/2}\,u\in X_{0}$. Since $X_{1}$ is
dense in $X_{0}$, the sequence $(v_{k})\subset X_{1}$ such that
$v_{k}\rightarrow v$ in $X_{0}$ as $k\rightarrow\infty$ exists.
From this and from (1.1) it follows that
$$
u_{k}:=(1+\psi^{2}(J))^{-1/2}\,v_{k}\rightarrow
u\quad\mbox{в}\quad X_{\psi}\quad\mbox{for}\quad
k\rightarrow\infty.
$$
It remains to note that
$$
u_{k}=(1+\psi^{2}(J))^{-1/2}J^{-1}Jv_{k}=J^{-1}(1+\psi^{2}(J))^{-1/2}Jv_{k}
\in X_{1}.
$$
Theorem 2.1 is proved.
\end{proof}

\begin{theorem}
Let functions $\psi,\chi\in\mathcal{B}$ be such that the function
$\psi/\chi$ is bounded in a neighborhood of $+\infty$. Then, for
each admissible couple $X=[X_{0},X_{1}]$ of Hilbert spaces, the
continuous and dense imbedding $X_{\chi}\hookrightarrow X_{\psi}$
holds. If the imbedding $X_{1}\hookrightarrow X_{0}$ is compact
and $\psi(t)/\chi(t)\rightarrow0$ as $t\rightarrow+\infty$, then
the imbedding $X_{\chi}\hookrightarrow X_{\psi}$ is also compact.
\end{theorem}

\begin{proof}
Let $J$ be a genarating operator for the couple $X$. Let us note
that $\mathrm{Spec}\,J\subseteq[r,+\infty)$ for some number $r>0$.
According to the condition of the theorem, we have
$\psi(t)/\chi(t)\leq c$ for $t\geq r$. Therefore
$$
X_{\chi}=\mathrm{Dom}\,\chi(J)\subseteq\mathrm{Dom}\,\psi(J)=X_{\psi},
\quad\|\psi(J)\,u\|_{X_{0}}\leq c\,\|\chi(J)\,u\|_{X_{0}}.
$$
From this formulae and from the definition of the spaces
$X_{\chi}$, $X_{\psi}$ we obtain the continuous imbedding
$X_{\chi}\hookrightarrow X_{\psi}$. Let us prove its density.

We consider the isometric isomorphisms
$\psi(J):X_{\psi}\leftrightarrow X_{0}$ and
$\chi(J):X_{\chi}\leftrightarrow X_{0}$. For any given $u\in
X_{\psi}$, we have $\psi(J)\,u\in X_{0}$. Since the space
$X_{\chi}$ is densely embedded into $X_{0}$, a sequence
$(v_{k})\subset X_{\chi}$ such that $v_{k}\rightarrow\psi(J)\,u$
in $X_{0}$ as $k\rightarrow\infty$ exists. Hence,
$\psi^{-1}(J)\,v_{k}\rightarrow u$ in $X_{\psi}$ as
$k\rightarrow\infty$, where
$$
\psi^{-1}(J)\,v_{k}=\psi^{-1}(J)\,\chi^{-1}(J)\,
\chi(J)\,v_{k}=\chi^{-1}(J)\,\psi^{-1}(J)\,\chi(J)\,v_{k} \in
X_{\chi}.
$$
Thus, we have proved the density of the embedding
$X_{\chi}\hookrightarrow X_{\psi}$.

Now let us assume that the imbedding $X_{1}\hookrightarrow X_{0}$
is compact and $\psi(t)/\chi(t)\rightarrow0$ as
$t\rightarrow+\infty$. We will prove the compactness of the
embedding $X_{\chi}\hookrightarrow X_{\psi}$. Let $(u_{k})$ be an
arbitrary bounded sequence belonging to $X_{\chi}$. Since the
sequence of elements $w_{k}:=J^{-1}\,\chi(J)\,u_{k}$ is bounded in
$X_{1}$, we can select a subsequence of elements
$w_{k_{n}}=J^{-1}\,\chi(J)\,u_{k_{n}}$ being the Cauchy sequence
in $X_{0}$. We show that $(u_{k_{n}})$ is the Cauchy sequence in
$X_{\psi}$.

Let $E_{t}$, $t\geq r$, be a resolution of the unity in $X_{0}$,
corresponding to the self-adjoint operator\;$J$. We can write
$$
\|u_{k_{n}}-u_{k_{m}}\|_{X_{\psi}}^{2}=
\|\psi(J)\,(u_{k_{n}}-u_{k_{m}})\|_{X_{0}}^{2}
=\|\psi(J)\,\chi^{-1}(J)\,J\,(w_{k_{n}}-w_{k_{m}})\|_{X_{0}}^{2}
$$
$$
=\int_{r}^{+\infty}\,\psi^{2}(t)\,\chi^{-2}(t)\,t^{2}\:
d\,\|E_{t}(w_{k_{n}}-w_{k_{m}})\|_{X_{0}}^{2}. \leqno(2.1)
$$
Let us choose an arbitrary number $\varepsilon>0$. There is a
number $\rho=\rho(\varepsilon)>r$ such that
$$
\psi(t)/\chi(t)\leq(2c_{0})^{-1}\varepsilon\quad\mbox{for}\quad
t\geq\rho\quad\mbox{and}\quad
c_{0}:=\sup\,\{\,\|w_{k}\|_{X_{1}}:\,k\in\mathbb{N}\,\}<\infty.
$$
Hence, for all indices $n,m$ we have
$$
\int_{\rho}^{+\infty}\psi^{2}(t)\,\chi^{-2}(t)\,t^{2}\:
d\,\|E_{t}(w_{k_{n}}-w_{k_{m}})\|_{X_{0}}^{2}\,\leq\,
(2c_{0})^{-2}\,\varepsilon^{2}\,\int_{\rho}^{+\infty}\,t^{2}\:
d\,\|E_{t}(w_{k_{n}}-w_{k_{m}})\|_{X_{0}}^{2}
$$
$$
\leq\,(2c_{0})^{-2}\,\varepsilon^{2}\,\|J\,(w_{k_{n}}-w_{k_{m}})\|_{X_{0}}^{2}\,
=\,(2c_{0})^{-2}\,\varepsilon^{2}\,\|w_{k_{n}}-w_{k_{m}}\|_{X_{1}}^{2}\,\leq
\,\varepsilon^{2} \leqno(2.2)
$$
In addition, by the inequality $\psi(t)/\chi(t)\leq c$ for $t\geq
r$, we can write the following:
$$
\int_{r}^{\rho}\,\psi^{2}(t)\,\chi^{-2}(t)\,t^{2}\:
d\,\|E_{t}(w_{k_{n}}-w_{k_{m}})\|_{X_{0}}^{2}\,\leq\,c^{2}\rho^{2}\,
\int_{r}^{\rho}\,d\,\|E_{t}(w_{k_{n}}-w_{k_{m}})\|_{X_{0}}^{2}
$$
$$
\leq
c^{2}\rho^{2}\,\|w_{k_{n}}-w_{k_{m}}\|_{X_{0}}^{2}\rightarrow0
\quad\mbox{as}\quad n,m\rightarrow\infty. \leqno(2.3)
$$
Now formulae (2.1) --- (2.3) imply the inequality
$\|u_{k_{n}}-u_{k_{m}}\|_{X_{\psi}}\leq2\varepsilon$ for
sufficiently large $n,m$. Therefore $(u_{k_{n}})$ is the Cauchy
sequence in the space $X_{\psi}$ which means the compactness of
the imbedding $X_{\chi}\hookrightarrow X_{\psi}$. Theorem 2.2 is
proved.
\end{proof}

\subsection{Reiteration.}

\begin{theorem}
Let functions $f,g,\psi\in\mathcal{B}$ be given. Suppose that the
function $f/g$ is bounded in a neighborhood of $+\infty$. Then
$[X_{f},X_{g}]_{\psi}=X_{\omega}$ holds with the equality of norms
for each admissible couple $X$ of Hilbert spaces. Here the
function $\omega\in\mathcal{B}$ is given by the formula
$\omega(t):=f(t)\,\psi(g(t)/f(t))$ for $t>0$. If $f,g,\psi$ are
interpolation parameters, so is $\omega$.
\end{theorem}

\begin{proof}
Since the function $f/g$ is bounded in a neighborhood of
$+\infty$, the couple $[X_{f},X_{g}]$ is admissible by Theorem 2.2
and, in addition, $\omega\in\mathcal{B}$. So, the spaces
$[X_{f},X_{g}]_{\psi}$ and $X_{\omega}$ are well defined. We will
prove them to be equal.

Let an operator $J$ be generating for the couple
$X=[X_{0},X_{1}]$, where $\mathrm{Spec}\,J\subseteq[r,+\infty)$
for some number $r>0$. We have three isometric isomorphisms
$$
f(J):X_{f}\leftrightarrow X_{0},\quad g(J):X_{g}\leftrightarrow
X_{0},\quad B:=f^{-1}(J)\,g(J):X_{g}\leftrightarrow X_{f}.
$$
Let us consider $B$ as a closed operator in the space $X_{f}$,
defined on $X_{g}$. The operator $B$ is generating for the couple
$[X_{f},X_{g}]$ because $B$ is positive and self-adjoint on
$X_{f}$. The positiveness of $B$ follows from the condition
$f(t)/g(t)\leq c$ for $t\geq r$ which implies
$$
(Bu,u)_{X_{f}}=(g(J)\,u,f(J)\,u)_{X_{0}}\geq
c^{-1}\,(f(J)\,u,f(J)\,u)_{X_{0}}= c^{-1}\,\|u\|_{X_{f}}^{2}.
$$
The self-adjointness follows from the fact that $0$ is a regular
point for the operator $B$.

Using the spectral theorem, we reduce the self-adjoint on $X_{0}$
operator $J$ to the form of multiplication by a function:
$J=I^{-1}(\alpha\cdot I)$. Here, $I:X_{0}\leftrightarrow
L_{2}(U,d\mu)$ is an isometric isomorphism, $(U,\mu)$ is a space
with a finite measure, $\alpha:U\rightarrow[r,+\infty)$ is a
measurable function. The isometric isomorphism
$If(J):X_{f}\leftrightarrow L_{2}(U,d\mu)$ reduces the
self-adjoint on $X_{f}$ operator $B$ to the form of multiplication
by the function $(g/f)\circ\alpha$:
$$
If(J)\,B\,u=Ig(J)\,u=(g\circ\alpha)\,Iu=(g\circ\alpha)\,If^{-1}(J)f(J)\,u=
((g/f)\circ\alpha)\,If(J)\,u,\quad u\in X_{g}.
$$
Therefore, for an arbitrary $u\in X_{\omega}$, we have
$$
\|\psi(B)\,u\|_{X_{f}}=
\|(\psi\circ(g/f)\circ\alpha)\cdot(If(J)\,u)\|_{L_{2}(U,d\mu)}=
\|(\omega\circ\alpha)\cdot(Iu)\|_{L_{2}(U,d\mu)}=\|\omega(J)\,u\|_{X_{0}}.
$$
Note that the function $f/\omega$ is bounded in a neighborhood of
$+\infty$. Hence, $X_{\omega}\hookrightarrow X_{f}$ and the
expression $f(J)\,u$ is well defined. Thus, the equality
$[X_{f},X_{g}]_{\psi}=X_{\omega}$ is proved.

We now assume that $f,g,\psi$ are interpolation parameters. We
show that $\omega$ is also an interpolation parameter. Let
arbitrary admissible couples $X=[X_{0},X_{1}]$, $Y=[Y_{0},Y_{1}]$
and a linear mapping $T$ be the same as those in Definition 2.2.
We have the bounded operators $T:X_{f}\rightarrow Y_{f}$ and
$T:X_{g}\rightarrow Y_{g}$ which imply the boundedness of the
operator $T:[X_{f},X_{g}]_{\psi}\rightarrow[Y_{f},Y_{g}]_{\psi}$.
We have already proved that $[X_{f},X_{g}]_{\psi}=X_{\omega}$ and
$[Y_{f},Y_{g}]_{\psi}=Y_{\omega}$. So, a bounded operator
$T:X_{\omega}\rightarrow Y_{\omega}$ exists. It means that
$\omega$ is an interpolation parameter. Theorem 2.3 is proved.
\end{proof}

\subsection{The interpolation of dual spaces.} Let $H$ be a
Hilbert space. We denote by $H'$ the space dual to $H$. Thus, $H'$
is the Banach space of all linear continuous functionals
$l:H\rightarrow\mathbb{C}$. By the Riesz theorem, the mapping
$S:v\mapsto(\,\cdot,v)_{H}$, where $v\in H$, establishes the
antilinear isometric isomorphism  $S:H\leftrightarrow H'$. This
implies that $H'$ is the Hilbert space with respect to the inner
product $(l,m)_{H'}:=(S^{-1}l,S^{-1}m)_{H}$. We emphasize that we
do not identify $H'$ as $H$ by means of the isomorphism $S$.

\begin{theorem}
Let $\psi\in\mathcal{B}$ be such that the function $\psi(t)/t$ is
bounded in a neighborhood of $+\infty$. Then, for each admissible
couple $[X_{0},X_{1}]$ of Hilbert spaces, the equality of spaces
$[X_{1}',X_{0}']_{\psi}=[X_{0},X_{1}]_{\chi}'$ with the equality
of norms holds. Here the function $\chi\in\mathcal{B}$ is given by
the formula $\chi(t):=t/\psi(t)$ for $t>0$. If $\psi$ is an
interpolation parameter, so is $\chi$.
\end{theorem}

\begin{proof}
Note that the couple $[X_{1}',X_{0}']$ is admissible, provided
that we naturally identify functionals from $X_{0}'$ as their
restrictions to the space $X_{1}$. From the condition of the
theorem it follows that $\varphi\in\mathcal{B}$. Thus, the spaces
$[X_{1}',X_{0}']_{\psi}$ and $[X_{0},X_{1}]_{\chi}'$ are well
defined. Let us prove these spaces to be equal.

Let $J:X_{1}\leftrightarrow X_{0}$ be a generating operator for
the couple $[X_{0},X_{1}]$. Let us consider the isometric
isomorphisms $S_{j}:X_{j}\leftrightarrow X_{j}'$, $j=0,\,1$, which
appear in the Riesz theorem. The operator $J'$, being adjoint to
$J$, satisfies the equality $J'=S_{1}J^{-1}S_{0}^{-1}$. This
results from the following:
$$
(J'l)u=l(Ju)=(Ju,S_{0}^{-1}l)_{X_{0}}=(u,J^{-1}S_{0}^{-1}l)_{X_{1}}=
(S_{1}J^{-1}S_{0}^{-1}l)u\quad\mbox{for each}\quad l\in X_{0}',\;
u\in X_{1}.
$$
Thus, the isometric isomorphism
$$
J'=S_{1}J^{-1}S_{0}^{-1}:\,X_{0}'\leftrightarrow X_{1}'
\leqno(2.4)
$$
exists.

Let us note that the equalities
$$
(u,JS_{1}^{-1}l)_{X_{0}}=(J^{-1}u,S_{1}^{-1}l)_{X_{1}}=l(J^{-1}u),\quad
(u,J^{-1}S_{0}^{-1}l)_{X_{0}}=(J^{-1}u,S_{0}^{-1}l)_{X_{0}}=l(J^{-1}u),
$$
where $l\in X_{0}'\hookrightarrow X_{1}'$, $u\in X_{0}$, imply the
property
$$
JS_{1}^{-1}l=J^{-1}S_{0}^{-1}l\in X_{1}\quad\mbox{for each}\quad
l\in X_{0}'.\leqno(2.5)
$$

Let us consider $J'$ as a closed operator in the space $X_{1}'$
with the domain $X_{0}'$. The operator $J'$ is generating for the
couple $[X_{1}',X_{0}']$ because $J'$ is positive and self-adjoint
on $X_{1}'$. The positiveness of $J'$ results from the
positiveness of the operator $J$ on the space $X_{0}$ and from
(2.5) in the following way:
$$
(J'l,l)_{X_{1}'}=(S_{1}J^{-1}S_{0}^{-1}l,l)_{X_{1}'}=
(J^{-1}S_{0}^{-1}l,S_{1}^{-1}l)_{X_{1}}=
(JJ^{-1}S_{0}^{-1}l,JS_{1}^{-1}l)_{X_{0}}
$$
$$
=(JJS_{1}^{-1}l,JS_{1}^{-1}l)_{X_{0}}\geq
c\,\|JS_{1}^{-1}l\|_{X_{0}}^{2}=c\,\|S_{1}^{-1}l\|_{X_{1}}^{2}=
c\,\|l\|_{X_{1}'}^{2}.
$$
Here the number $c>0$ does not depend on $l\in X_{0}'$. The
operator $J'$ is self-adjoint because $0$ is its regular point
(see (2.4)). Let us reduce the operator $J$ to the form of
multiplication by a function: $J=I^{-1}(\alpha\cdot I)$ as it has
been done in the proof of Theorem~2.3. The isometric isomorphism
$$
IJS_{1}^{-1}:\,X_{1}'\leftrightarrow L_{2}(U,d\mu) \leqno(2.6)
$$
reduces the operator $J'$ to the form of multiplication by the
same function $\alpha$:
$$
(IJS_{1}^{-1})J'l=IS_{0}^{-1}l=IJJ^{-1}S_{0}^{-1}l=\alpha\cdot
IJ^{-1}S_{0}^{-1}l=\alpha\cdot IJS_{1}^{-1}l\quad\mbox{for
each}\quad l\in X_{0}'.
$$
The last equality follows from (2.5).

By Theorem 2.2, two continuous dense imbeddings
$X_{0}'\hookrightarrow[X_{1}',X_{0}']_{\psi}$ and
$[X_{0},X_{1}]_{\chi}\hookrightarrow X_{0}$ hold. The second
imbedding implies the continuous dense imbedding
$X_{0}'\hookrightarrow[X_{0},X_{1}]_{\chi}'$. Let us show that the
norms in the spaces $[X_{1}',X_{0}']_{\psi}$ and
$[X_{0},X_{1}]_{\chi}'$ are equal on the dense subset $X_{0}'$.
For each $l\in X_{0}'$, $u\in[X_{0},X_{1}]_{\chi}$, we can write
$$
l(u)=(u,S_{0}^{-1}l)_{X_{0}}=(\chi(J)u,\chi^{-1}(J)S_{0}^{-1}l)_{X_{0}}=
(v,\chi^{-1}(J)S_{0}^{-1}l)_{X_{0}}
$$
with $v:=\chi(J)u\in X_{0}$. It implies the following:
$$
\|l\|_{\,[X_{0},X_{1}]_{\chi}'}=
\sup\,\{\,|l(u)|\,/\,\|u\|_{\,[X_{0},X_{1}]_{\chi}}: u\in
[X_{0},X_{1}]_{\chi},\,u\neq0\,\}
$$
$$
=\sup\,\{\,|(v,\chi^{-1}(J)S_{0}^{-1}l)_{X_{0}}|\,/\,\|v\|_{X_{0}}:
\,v\in X_{0},\,v\neq0\,\}
$$
$$
=\|\chi^{-1}(J)S_{0}^{-1}l\|_{X_{0}}=
\|I\chi^{-1}(J)S_{0}^{-1}l\|_{L_{2}(U,d\mu)}=
\|(\chi^{-1}\circ\alpha)\cdot IS_{0}^{-1}l\|_{L_{2}(U,d\mu)}.
$$
On the other hand, using isomorphisms (2.6), (2.4), we have
$$
\|l\|_{\,[X_{1}',X_{0}']_{\psi}}=\|\psi(J')l\|_{X_{1}'}=
\|\chi^{-1}(J')J'l\|_{X_{1}'}=
\|(IJS_{1}^{-1})\chi^{-1}(J')J'l\|_{L_{2}(U,d\mu)}
$$
$$
=\|(\chi^{-1}\circ\alpha)\cdot(IJS_{1}^{-1})J'l\|_{L_{2}(U,d\mu)}
=\|(\chi^{-1}\circ\alpha)\cdot IS_{0}^{-1}l\|_{L_{2}(U,d\mu)}.
$$
Thus, norms in the spaces $[X_{1}',X_{0}']_{\psi}$ and
$[X_{0},X_{1}]_{\chi}'$ are equal on the dense subset $X_{0}'$.
So, these spaces coincide.

Now suppose $\psi$ to be an interpolation parameter. We will show
that so is $\chi$. Let admissible couples $X=[X_{0},X_{1}]$,
$Y=[Y_{0},Y_{1}]$ and a linear mapping $T$ be the same as those in
Definition 2.2. Passing to the adjoint operator $T'$, we get the
bounded operators $T':Y_{j}'\rightarrow X_{j}'$, $j=0,\,1$. Since
$\psi$ is an interpolation parameter, a bounded operator
$T':[Y_{1}',Y_{0}']_{\psi}\rightarrow[X_{1}',X_{0}']_{\psi}$
exists. As we have already proved,
$[X_{1}',X_{0}']_{\psi}=[X_{0},X_{1}]_{\chi}'$ and
$[Y_{1}',Y_{0}']_{\psi}=[Y_{0},Y_{1}]_{\chi}'$ with equalities of
norms. Therefore a bounded operator
$T':[Y_{0},Y_{1}]_{\chi}'\rightarrow[X_{0},X_{1}]_{\chi}'$ exists.
Thus, passing to the second adjoint operator $T''$, we get the
bounded operator
$T'':[X_{0},X_{1}]_{\chi}''\rightarrow[Y_{0},Y_{1}]_{\chi}''$. It
remains to identify the second dual spaces with original spaces
which leads us to the bounded operator
$T:[X_{0},X_{1}]_{\chi}\rightarrow[Y_{0},Y_{1}]_{\chi}$. This
means that $\chi$ is an interpolation parameter. Theorem 2.4 is
proved.
\end{proof}

\subsection{The ihterpolation of direct products of
spaces.}

\begin{theorem}
Let a finite or countable set of admissible couples of Hilbert
spaces $X^{(k)}:=[X_{0}^{(k)},X_{1}^{(k)}]$, $k\in\omega$, be
given. Suppose that the set of norms of the imbedding operators
$X_{1}^{(k)}\hookrightarrow X_{0}^{(k)}$, $k\in\omega$, is
bounded. Then, for an arbitrary function $\psi\in\mathcal{B}$, the
equality of spaces
$$
\left[\:\prod_{k\in\omega}X_{0}^{(k)},\,
\prod_{k\in\omega}X_{1}^{(k)}\right]_{\psi}=
\,\prod_{k\in\omega}\left[X_{0}^{(k)},X_{1}^{(k)}\right]_{\psi}
$$
and the equality of norms in them hold.
\end{theorem}

\begin{proof}
We assume that $\omega=\mathbb{N}$ (the
case of finite set $\omega$ is treated analogously and easier).
The spaces $X_{0}:=\prod_{k=1}^{\infty}X_{0}^{(k)}$,
$X_{1}:=\prod_{k=1}^{\infty}X_{1}^{(k)}$ are Hilbert and separable
ones. The continuous imbedding $X_{1}\hookrightarrow X_{0}$ is
evident due to the condition of the theorem. Let
$u:=(u_{1},u_{2},\ldots)\in X_{0}$. For all indices $n,k$ an
element $v_{n,k}\in X_{1}^{(k)}$ such that
$\|u_{k}-v_{n,k}\|_{X_{0}^{(k)}}<1/n$ exists. Let us form a
sequence of vectors
$v^{(n)}:=(v_{n,1},\ldots,v_{n,n},0,0,\ldots)\in X_{1}$. We have
$$
\|u-v^{(n)}\|_{X_{0}}^{2}=\sum_{k=1}^{n}\|u_{k}-v_{n,k}\|_{X_{0}}^{2}+
\sum_{k=n+1}^{\infty}\|u_{k}\|_{X_{0}}^{2}\leq \frac{n}{n^{2}}+
\sum_{k=n+1}^{\infty}\|u_{k}\|_{X_{0}}^{2}\rightarrow0\;\;\mbox{as}
\;\;n\rightarrow\infty.
$$
Thus, the couple $X:=[X_{0},X_{1}]$ is admissible.

Denote by $J_{k}$ a generating operator for the couple $X^{(k)}$.
An operator $J:=(J_{1},J_{2},\ldots)$ is generating for the couple
$X$ which may be proved directly. Moreover, it is natural to
expect that $\psi(J)=(\psi(J_{1}),\psi(J_{2}),\ldots)$ and
$\mathrm{Dom}\,\psi(J)=\prod_{k=1}^{\infty}X_{\psi}^{(k)}$. Now we
prove these equalities. Let us reduce the operator $J_{k}$ to the
form of multiplication by a function: $I_{k}J_{k}=\alpha_{k}\cdot
I_{k}$. Here $I_{k}:X_{0}^{(k)}\leftrightarrow
L_{2}(V_{k},d\mu_{k})$ is an isometric isomorphism, $V_{k}$ is a
space with a finite measure $\mu_{k}$ and
$\alpha_{k}:V_{k}\rightarrow(0,+\infty)$ is a measurable function.
We may consider the sets $V_{k}$ to be mutually disjoint. Let us
set $V:=\bigcup_{k=1}^{\infty}V_{k}$. We call $\Omega\subseteq V$
a measurable set if, for every index $k$, the set $\Omega\cap
V_{k}$ is $\mu_{k}$-measurable. On the $\sigma$-algebra of all
measurable sets $\Omega\subseteq V$, we introduce the
$\sigma$-finite measure
$\mu(\Omega):=\sum_{k=1}^{\infty}\mu_{k}(\Omega\cap V_{k})$.
Further, for a vector $u:=(u_{1},u_{2},\ldots)\in X_{0}$, we
consider the measurable functions $Iu$ and $\alpha$, defined on
the set $V$ by the formulae $(Iu)(\lambda):=(I_{k}u_{k})(\lambda)$
and $\alpha(\lambda):=\alpha_{k}(\lambda)$ with $\lambda\in
V_{k}$. Now we have the isometric isomorphism
$I:X_{0}\leftrightarrow L_{2}(V,d\mu)$. It reduces the operator
$J$ to the form of multiplication by the function $\alpha$ because
$$
(IJu)(\lambda)=(I_{k}J_{k}u_{k})(\lambda)=
\alpha_{k}(\lambda)(I_{k}u_{k})(\lambda)=\alpha(\lambda)(Iu)(\lambda)\quad
\mbox{for}\;\;u\in X_{1},\,\lambda\in V_{k}.
$$
Hence, we can write down the following:
$$
X_{\psi}=\mathrm{Dom}\,\psi(J)=\{\,u\in
X_{0}:\,(\psi\circ\alpha)\cdot(Iu)\in L_{2}(V,d\mu)\,\}
$$
$$
=\{\,u\in
X_{0}:\,\sum_{k=1}^{\infty}\:\|(\psi\circ\alpha_{k})\cdot(I_{k}u_{k})\|
_{L_{2}(V_{k},d\mu_{k})}^{2}<\infty\,\}
$$
$$
=\{\,u:\,u_{k}\in\mathrm{Dom}\,\psi(J_{k}),\;
\sum_{k=1}^{\infty}\:\|\psi(J_{k})u_{k}\|
_{X_{0}^{(k)}}^{2}<\infty\,\}=\prod_{k=1}^{\infty}X_{\psi}^{(k)}.
$$
Furthermore, for each $u\in \mathrm{Dom}\,\psi(J)$, we have
$$
(I\psi(J)u)(\lambda)=\psi(\alpha(\lambda))\,(Iu)(\lambda)=
\psi(\alpha_{k}(\lambda))\,(I_{k}u_{k})(\lambda)
$$
$$
=(I_{k}\psi(J_{k})u_{k})(\lambda)=
\bigl(I(\psi(J_{1})u_{1},\psi(J_{2})u_{2},\ldots)\bigr)(\lambda)
\quad\mbox{for}\;\;\lambda\in V_{k}.
$$
Therefore $\psi(J)u=(\psi(J_{1})u_{1},\psi(J_{2})u_{2},\ldots)$
which implies

$$
\|u\|_{X_{\psi}}^{2}=\|\psi(J)u\|_{X_{0}}^{2}=
\sum_{k=1}^{\infty}\,\|\psi(J_{k})u_{k}\|_{X_{0}^{(k)}}^{2}=
\sum_{k=1}^{\infty}\,\|u_{k}\|_{X_{\psi}^{(k)}}^{2}.
$$
Theorem 2.5 is proved.
\end{proof}

\subsection{An operator norm in interpolation spaces.}

\begin{theorem}
For given interpolation parameter $\psi\in\mathcal{B}$ and number
$m>0$, there exists a number $c=c(\psi,m)>0$ such that
$$
\|T\|_{X_{\psi}\rightarrow Y_{\psi}}\leq c
\max\,\bigl\{\,\|T\|_{X_{j}\rightarrow Y_{j}}:\,j=0,\,1\,\bigr\}.
$$
Here $X=[X_{0},X_{1}]$ and $Y=[Y_{0},Y_{1}]$ are admissible
couples of Hilbert spaces for which the norms of the imbedding
operators $X_{1}\hookrightarrow X_{0}$ and $Y_{1}\hookrightarrow
Y_{0}$ do not exceed the number $m$, and $T$ is any linear mapping
defined on the space $X_{0}$ and establishing the bounded
operators $T:X_{j}\rightarrow Y_{j}$ with $j=0,\,1$.
\end{theorem}

\begin{proof}
Let us suppose the contrary. Then we can write:
$$
\|T_{k}\|_{X^{(k)}_{\psi}\rightarrow Y^{(k)}_{\psi}}>k\,m_{k}
\quad\mbox{for each index}\;\;k. \leqno(2.7)
$$
Here, $X^{(k)}:=[X_{0}^{(k)},X_{1}^{(k)}]$ and
$Y^{(k)}:=[Y_{0}^{(k)},Y_{1}^{(k)}]$ are some admissible couples
of Hilbert spaces for which the norms of the imbedding operators
$X_{1}^{(k)}\hookrightarrow X_{0}^{(k)}$ and
$Y_{1}^{(k)}\hookrightarrow Y_{0}^{(k)}$ do not exceed the number
$m$. Furthermore, $T_{k}$ is a certain linear mapping defined on
the space $X_{0}^{(k)}$ and establishing the bounded operators
$T_{k}:X_{j}^{(k)}\rightarrow Y_{j}^{(k)}$ with $j=0,\,1$. We also
use the notation
$$
m_{k}:=\max\left\{\,\|T_{k}\|_{X^{(k)}_{0}\rightarrow
Y^{(k)}_{0}},\;\|T_{k}\|_{X^{(k)}_{1}\rightarrow
Y^{(k)}_{1}}\,\right\}>0.
$$
Now let us consider the bounded operators
$$
T:u=(u_{1},u_{2},\ldots)\mapsto
(m_{1}^{-1}\,T_{1}u_{1},m_{2}^{-1}\,T_{2}u_{2},\ldots),\quad
T:\,\prod_{k=1}^{\infty}X_{j}^{(k)}\rightarrow\prod_{k=1}^{\infty}Y_{j}^{(k)},
\;j=0,\,1. \leqno(2.8)
$$
Their boundedness results from the following inequalities:
$$
\sum_{k=1}^{\infty}\;\left\|m_{k}^{-1}\,T_{k}u_{k}\right\|_{Y_{j}^{(k)}}^{2}\leq
\sum_{k=1}^{\infty}\;m_{k}^{-2}\,\|T_{k}\|_{X_{j}^{(k)}\rightarrow
Y_{j}^{(k)}}^{2}\,\|u_{k}\|_{X_{j}^{(k)}}^{2}\leq
\sum_{k=1}^{\infty}\;\|u_{k}\|_{X_{j}^{(k)}}^{2}.
$$
Since $\psi$ is an interpolation parameter, the boundedness of
operators (2.8) implies the existence of the bounded operator
$$
T:\,\left[\:\prod_{k=1}^{\infty}X_{0}^{(k)},\,
\prod_{k=1}^{\infty}X_{1}^{(k)}\right]_{\psi}\rightarrow\;
\left[\:\prod_{k=1}^{\infty}Y_{0}^{(k)},\,
\prod_{k=1}^{\infty}Y_{1}^{(k)}\right]_{\psi}
$$
which by Theorem 2.5 means the boundedness of the operator
$$
T:\,\prod_{k=1}^{\infty}X_{\psi}^{(k)}\rightarrow\;
\prod_{k=1}^{\infty}Y_{\psi}^{(k)}.
$$
Let $c_{0}$ be the norm of the last operator. For every index $k$
we consider a vector $u^{(k)}:=(u_{1},\ldots,u_{k},\ldots)$ such
that $u_{k}\in X_{\psi}^{(k)}$ and $u_{j}=0$ for $j\neq k$. We
have:
$$
\|T_{k}u_{k}\|_{Y_{\psi}^{(k)}}=
m_{k}\,\|Tu^{(k)}\|_{\prod_{j=1}^{\infty}Y_{\psi}^{(j)}}\leq
m_{k}\,c_{0}\,\|u^{(k)}\|_{\prod_{j=1}^{\infty}X_{\psi}^{(j)}}=
m_{k}\,c_{0}\,\|u_{k}\|_{X_{\psi}^{(k)}}
$$
for each $u_{k}\in X_{\psi}^{(k)}$. Hence,
$$
\|T_{k}\|_{X^{(k)}_{\psi}\rightarrow Y^{(k)}_{\psi}}\leq
c_{0}\,m_{k} \quad\mbox{for every index}\;\;k,
$$
contrary to inequality (2.7). Thus, our supposition is false and
the theorem is true.
\end{proof}

\subsection{A criterion for a function to be an
interpolation parameter.} Using Peetre's results \cite{Peet68},
\cite[Sec. 5.4]{BL80}, we prove the following criterion.

\begin{definition}
Let a function
$\psi:(0,+\infty)\rightarrow(0,+\infty)$ and a number $r\geq0$ be
given. The function $\psi$ is called \textit{quasiconcave} (or
pseudoconcave) on the semiaxis $(r,+\infty)$ if a concave function
$\psi_{1}:(r,+\infty)\rightarrow(0,+\infty)$ such that
$\psi(t)\asymp \psi_{1}(t)$ for $t>r$ exists. The function $\psi$
is called quasiconcave in a neighborhood of $+\infty$ if it is
quasiconcave on a certain semiaxis $(r,+\infty)$, where $r$ is a
sufficiently large number.
\end{definition}

\begin{theorem}
A function $\psi\in\mathcal{B}$
is an interpolation parameter if and only if it is quasiconcave in
a neighborhood of $+\infty$.
\end{theorem}

To prove this theorem we need two lemmas.

\begin{lemma}
Let a function $\psi$ belong to the set $\mathcal{B}$ and be
quasiconcave in a neighborhood of $+\infty$. Then there exists a
concave function $\psi_{0}:(0,+\infty)\rightarrow(0,+\infty)$ such
that for every number $\varepsilon>0$ it holds $\psi(t)\asymp
\psi_{0}(t)$ with $t\geq\varepsilon$.
\end{lemma}

\begin{proof} It is evident.
\end{proof}

\begin{lemma}
Let a function $\psi\in\mathcal{B}$ and a number $r\geq0$ be
given. The function $\psi$ is quasiconcave on the semiaxis
$(r,+\infty)$ if and only if there exists a number $c>0$ such that
$$
\psi(t)/\psi(s)\leq c\,\max\{1,\,t/s\}\quad\mbox{for each}\quad
t,s>r.
$$
\end{lemma}

\begin{proof}
In the case where $r=0$ this lemma was proved by J. Peetre
\cite{Peet68}, \cite[Lemma 5.4.3]{BL80} (the condition
$\psi\in\mathcal{B}$ being superfluous). In the case where $r>0$
the sufficiency can be proved analogously. The necessity is be
reduced to the case $r=0$ with the help of Lemma 2.1. Indeed, let
us put $\varepsilon=r$ in this lemma. Then we have a function
$\psi_{0}$ such that
$$
\psi(t)/\psi(s)\asymp\psi_{0}(t)/\psi_{0}(s)\leq
c_{0}\max\{1,\,t/s\}\quad\mbox{for each}\quad t,s>r.
$$
(In fact, $c_{0}=1$ for a concave function $\psi_{0}$
\cite{Peet68}). Lemma 2.2 is proved.
\end{proof}

\begin{proof}[Proof of theorem $2.7$]
\textit{Sufficiency.} Let us suppose that a function
$\psi\in\mathcal{B}$ is quasiconcave in a neighborhood of
$+\infty$. We need to prove that $\psi$ is an interpolation
parameter.

Let admissible couples $X=[X_{0},X_{1}]$, $Y=[Y_{0},Y_{1}]$ and a
linear mapping $T$ be the same as those in Definition 2.2. In
addition, let operators $J_{X}:X_{1}\leftrightarrow X_{0}$ and
$J_{Y}:Y_{1}\leftrightarrow Y_{0}$ be generating ones for the
couples $X$ and $Y$ respectively. Using the spectral theorem we
reduce these operators, self-adjoint in $X_{0}$ and in $Y_{0}$
respectively, to the form of multiplication by a function:
$$
J_{X}=I_{X}^{-1}\,(\alpha\cdot I_{X})\quad\mbox{and}\quad
J_{Y}=I_{Y}^{-1}\,(\beta\cdot I_{Y}). \leqno(2.9)
$$
Here, $I_{X}:X_{0}\leftrightarrow L_{2}(U,d\mu)$ and
$I_{Y}:Y_{0}\leftrightarrow L_{2}(V,d\nu)$ are certain isometric
isomorphisms, $(U,\mu)$ and $(V,\nu)$ are spaces with finite
measures and $\alpha:U\rightarrow(0,+\infty)$ and
$\beta:V\rightarrow(0,+\infty)$ are some measurable functions.
Since the operators $T:X_{0}\rightarrow Y_{0}$ and
$T:X_{1}\rightarrow Y_{1}$ are bounded, so are the operators
$$
I_{Y}\,T\,I_{X}^{-1}:\,L_{2}(U,d\mu)\rightarrow L_{2}(V,d\nu),
\leqno(2.10)
$$
$$
I_{Y}\,J_{Y}\,T\,J_{X}^{-1}\,I_{X}^{-1}:\,L_{2}(U,d\mu)\rightarrow
L_{2}(V,d\nu). \leqno(2.11)
$$
By virtue of (2.9), we can write
$$
I_{Y}\,J_{Y}\,T\,J_{X}^{-1}\,I_{X}^{-1}=(\beta\cdot
I_{Y})\,T\,(\alpha^{-1}\cdot I_{X}^{-1}).
$$
Hence (2.11) implies the boundedness of the operator
$$
I_{Y}\,T\,I_{X}^{-1}=\beta^{-1}\cdot(I_{Y}\,J_{Y}\,T\,J_{X}^{-1}\,I_{X}^{-1})
\cdot\alpha:\,L_{2}(U,\alpha^{2}d\mu)\rightarrow
L_{2}(V,\beta^{2}d\nu). \leqno(2.12)
$$

Let a concave function
$\psi_{0}:(0,+\infty)\rightarrow(0,+\infty)$ be the same as that
in Lemma 2.1. Note that  $\psi_{0}\in\mathcal{B}$ and (see
Subsection 2.1)
$$
X_{\psi}=X_{\psi_{0}},\quad Y_{\psi}=Y_{\psi_{0}}\quad\mbox{with
equivalence of norms}. \leqno(2.13)
$$
J. Peetre \cite{Peet68}, \cite[Theorem 5.4.4]{BL80} proved that a
positive function is quasiconcave on $(0,+\infty)$ if and only if
it is an interpolation function in the sense of the definition
stated in \cite[Definition 5.4.2]{BL80}. Hence, for the function
$\psi_{0}$, the boundedness of operators (2.10), (2.12) implies
the existence of a bounded operator
$$
I_{Y}\,T\,I_{X}^{-1}:\,L_{2}(U,(\psi_{0}\circ\alpha^{2})\,d\mu)\rightarrow
L_{2}(V,(\psi_{0}\circ\beta^{2})\,d\nu). \leqno(2.14)
$$
Let us pass from (2.14) to the operator $T:X_{\psi_{0}}\rightarrow
Y_{\psi_{0}}$ with the help of the isometric isomorphisms
$\psi_{0}(J_{X}):\,X_{\psi_{0}}\leftrightarrow X_{0}$ and
$\psi_{0}(J_{Y}):\,Y_{\psi_{0}}\leftrightarrow Y_{0}$. We reduce
these isomorphisms (which are self-adjoint operators on $X_{0}$
and $Y_{0}$ respectively) to the form of the multiplication by a
function:
$$
I_{X}\,\psi_{0}(J_{X})=(\psi_{0}\circ\alpha)\cdot I_{X}:\,
X_{\psi_{0}}\leftrightarrow L_{2}(U,d\mu),
$$
$$
I_{Y}\,\psi_{0}(J_{Y})=(\psi_{0}\circ\beta)\cdot I_{Y}:\,
Y_{\psi_{0}}\leftrightarrow L_{2}(V,d\nu).
$$
We get the isometric isomorphisms
$$
I_{X}=(\psi_{0}^{-1}\circ\alpha)\cdot(I_{X}\,\psi_{0}(J_{X})):\,
X_{\psi_{0}}\leftrightarrow L_{2}(U,(\psi^{2}\circ\alpha)\,d\mu),
$$
$$
I_{Y}=(\psi_{0}^{-1}\circ\beta)\cdot(I_{Y}\,\psi_{0}(J_{Y})):\,
Y_{\psi_{0}}\leftrightarrow L_{2}(V,(\psi^{2}\circ\beta)\,d\nu).
$$
From this and (2.14) the existence of the bounded operator
$$
T=I_{Y}^{-1}(I_{Y}\,T\,I_{X}^{-1})I_{X}:\,X_{\psi_{0}}\rightarrow
Y_{\psi_{0}}
$$
follows.

Thus, due to equations (2.13) we have
$$
(T:X_{j}\rightarrow
Y_{j},\,j=0,1)\,\Rightarrow\,(T:X_{\psi_{0}}\rightarrow
Y_{\psi_{0}})\,\Rightarrow\,(T:X_{\psi}\rightarrow Y_{\psi}),
$$
where the linear operators are bounded. So, by Definition 2.2 the
function $\psi$ is an interpolation parameter. Sufficiency is
proved.

\textit{Necessity.} Now we suppose that a function
$\psi\in\mathcal{B}$ is an interpolation parameter. We need to
prove that $\psi$ is quasiconcave in a neighborhood of $+\infty$.
The proof is similar to \cite{Peet68}, \cite[Sec. 5.4]{BL80}.

Let us consider the space $L_{2}(U,d\mu)$ with $U=\{0,\,1\}$,
$\mu(\{0\})=\mu(\{1\})=1$ and define on it the linear mapping $T$
by the formula $(Tu)(0):=0$, $(Tu)(1):=u(0)$, where $u\in
L_{2}(U,d\mu)$. Choose arbitrary numbers $s,t>1$ and put
$\omega(0):=s^{2}$, $\omega(1):=t^{2}$. We have the admissible
couple $X:=[L_{2}(U,d\mu),L_{2}(U,\omega\,d\mu)]$ and bounded
operators
$$
T:\,L_{2}(U,d\mu)\rightarrow L_{2}(U,d\mu)\quad\mbox{and}\quad
T:\,L_{2}(U,\omega\,d\mu)\rightarrow L_{2}(U,\omega\,d\mu)
$$
with norms $1$ and $t/s$ respectively. From this, since $\psi$ is
an interpolation parameter, it follows that a bounded operator
$T:X_{\psi}\rightarrow X_{\psi}$ exists. By Theorem\;2.6 with
$Y=X$ and $m=1$ we conclude that the norm of this operator
satisfies the inequality
$$
\|T\|_{X_{\psi}\rightarrow X_{\psi}}\leq\,c\,\max\{1,t/s\}.
\leqno(2.15)
$$
Here, the number $c>0$ does not depend on $t,s>1$.

It is not difficult to calculate the norm in the space $X_{\psi}$.
Indeed, the operator $J$ of multiplication by the function
$\omega^{1/2}$ is generating for the couple $X$. Hence, since
$\psi(J)$ is the operator of multiplication by the function
$\psi\circ\omega^{1/2}$, we can write
$$
\bigl\|u\bigr\|_{X_{\psi}}^{2}=
\bigl\|(\psi\circ\omega^{1/2})\,u\bigr\|_{L_{2}(U,d\mu)}^{2}=
\psi^{2}(s)\,|u(0)|^{2}+\psi^{2}(t)\,|u(1)|^{2},\quad
\|Tu\|_{X_{\psi}}^{2}=\psi^{2}(t)\,|u(0)|^{2}.
$$
It follows that
$$
\|T\|_{X_{\psi}\rightarrow X_{\psi}}=\psi(t)/\psi(s). \leqno(2.16)
$$
Now relations (2.15), (2.16) imply the inequality
$$
\psi(t)\leq c\,\max\{1,t/s\}\,\psi(s)\quad\mbox{for each}\quad
t,s>1.
$$
According to Theorem 2.2, the last statement is equivalent to the
quasiconcavity of the function $\psi$ on the semiaxis
$(1,+\infty)$. Necessity is proved.
\end{proof}

\section{A refined scale of spaces}

\subsection{Quasiregularly varying functions.} We recall
the following:

\begin{definition}
A positive function $\psi$ defined on a semiaxis $[b,+\infty)$ is
called a function \textit{regularly varying} at $+\infty$ with the
index $\theta\in\mathbb{R}$ if $\psi$ is Borel measurable on
$[b_{\,0},+\infty)$ for some number $b_{\,0}\geq b$ and
$$
\lim_{t\rightarrow\,+\infty}\;\psi(\lambda\,t)/\psi(t)=
\lambda^{\theta}\quad\mbox{for each}\quad \lambda>0.
$$
A function regularly varying at $+\infty$ with the index
$\theta=0$ is called a function \textit{slowly varying} at
$+\infty$.
\end{definition}

The theory of regularly varying functions was founded by J.
Karamata in the 1930s. These functions are closely related to the
power functions and have numerous applications, mainly due to
their special role in Tauberian-type theorems \cite{Se85, Re87,
BiGoTe89, Ma00}. A standard example of functions regularly varying
at $+\infty$ with the index $\theta$ is
$$
\psi(t)=t^{\theta}\,(\ln t)^{r_{1}}\,(\ln\ln t)^{r_{2}} \ldots
(\ln\ldots\ln t)^{r_{k}}\quad\mbox{for}\quad t\gg1,
$$
where $r_{1}, r_{2},\ldots,r_{k}\in\mathbb{R}$. In the case where
$\theta=0$ these functions form the \textit{logarithmic
multiscale} which has a number of applications in the theory of
function spaces.

\begin{definition}
A positive function $\psi$ defined on a semiaxis $[b,+\infty)$ is
called a function \textit{quasiregularly varying} at $+\infty$
with the index $\theta\in\mathbb{R}$ if there exist a number
$b_{1}\geq b$ and a function $\psi_{1}:[b_{1},+\infty)\rightarrow
(0,+\infty)$ regularly varying at $+\infty$ with the index
$\theta$ such that $\psi(t)\asymp\psi_{1}(t)$ with $t\geq b_{1}$.
A function quasiregularly varying  at $+\infty$ with the index
$\theta=0$ is called a function \textit{quasislowly varying} at
$+\infty$.
\end{definition}

We denote by $\mathrm{QSV}$ the set of all functions quasislowly
varying at $+\infty$. It is evident that $\psi$ is a function
quasiregularly varying at $+\infty$ with the index $\theta$ if and
only if $\psi(t)=t^{\theta}\varphi(t)$, $t\gg1$, for some function
$\varphi\in\mathrm{QSV}$. From the known \cite[Theorem 1.2]{Se85}
integral representation of a slowly varying function it
immediately results the following description of the set
$\mathrm{QSV}$.

\begin{theorem}
Let $\varphi\in\mathrm{QSV}$.
Then
$$
\varphi(t)=
\exp\left(\beta(t)+\int_{r}^{\:t}\frac{\alpha(\tau)}{\tau}\;d\tau\right),
\quad t\geq r, \leqno(3.1)
$$
for some number $r>0$, a continuous function $\alpha:
[r,+\infty)\rightarrow\mathbb{R}$ approaching zero at $+\infty$
and a bounded function $\beta: [r,+\infty)\rightarrow\mathbb{R}$.
The converse statement is also true: every function of form \rm
(3.1) \it belongs to the set $\mathrm{QSV}$.
\end{theorem}

Following interpolation property of quasiregular varying functions
will play a decisive role in further.

\begin{theorem}
Let a function $\psi\in\mathcal{B}$ be quasiregularly varying at
$+\infty$ with the index $\theta\in(0,1)$. Then $\psi$ is an
interpolation parameter.
\end{theorem}

\begin{proof}
We can write $\psi(t)=t^{\theta}\varphi(t)$ for $t>0$ with
$\varphi\in\mathrm{QSV}$. According to Theorem 3.1, the function
$\varphi$ can be represented in form (3.1). Let us set
$\varepsilon:=\min\{\theta,1-\theta\}>0$ and choose a number
$r_{\varepsilon}\geq r$ such that $|\alpha(t)|<\varepsilon$ for
$t> r_{\varepsilon}$. For each $t,s>r_{\varepsilon}$, we have by
virtue of (3.1) the following:
$$
\frac{\varphi(t)}{\varphi(s)}=\exp\left(\beta(t)-\beta(s)+
\int_{s}^{\:t}\frac{\alpha(\tau)}{\tau}\;d\tau\right) \leq
c\exp\left|\int_{s}^{\:t}\frac{\varepsilon}{\tau}\;d\tau\right|=
c\max\left\{(t/s)^{\varepsilon},(s/t)^{\varepsilon}\right\}.
$$
Here the number $c>0$ does not depend on $t$ and $s$ because the
function $\beta$ is bounded. From this and from the inequality
$0\leq\theta\pm\varepsilon\leq1$ it follows that
$$
\psi(t)/\psi(s)=(t^{\theta}\varphi(t))/(s^{\theta}\varphi(s))\leq
c\max\left\{(t/s)^{\theta+\varepsilon},(t/s)^{\theta-\varepsilon}\right\}\leq
c\max\{1,t/s\}.
$$
Hence, by Theorem 2.2 the function $\psi\in\mathcal{B}$ is
quasiconcave in a neighborhood of $+\infty$. According to Theorem
2.7, this is equivalent to the statement that $\psi$ is an
interpolation parameter. Theorem 3.2 is proved.
\end{proof}

We need the following properties of the set $\mathrm{QSV}$.

\begin{theorem}
Let $\varphi,\chi\in\mathrm{QSV}$. The following assertions are
true.
\begin{itemize}
  \item [(i)] \it There is a positive function $\varphi_{1}\in
C^{\infty}((0;+\infty))$ regularly varying at $+\infty$ such that
$\varphi(t)\asymp\varphi_{1}(t)$ with $t\gg1$.
  \item [(ii)] \it For each $\theta>0$, the limits
$t^{-\theta}\varphi(t)\rightarrow0$ and
$t^{\theta}\varphi(t)\rightarrow+\infty$ as $t\rightarrow+\infty$
hold.
  \item [(iii)] \it The functions $\varphi+\chi$, $\varphi\,\chi$, $\varphi/\chi$
and $\varphi^{\sigma}$, where $\sigma\in\mathbb{R}$, belong to the
set $\mathrm{QSV}$.
  \item [(iv)] \it Let $\theta\geq0$ and in the case where $\theta=0$ suppose
that $\varphi(t)\rightarrow\infty$ as $t\rightarrow+\infty$. Then
the composite function $\chi(t^{\theta}\varphi(t))$ of the
argument $t$ belongs to the set $\mathrm{QSV}$.
  \end{itemize}
\end{theorem}

\begin{proof}
For regularly varying functions
$\varphi,\chi$ these assertions are known \cite[Sec. 1.5]{Se85}
(even with the strong equivalence being in assertion (i)). This
implies immediately assertions (i), (ii), (iii) for the functions
$\varphi,\chi\in\mathrm{QSV}$.

It remains to prove assertion (iv). Let $\lambda>0$. Since
$\varphi\in\mathrm{QSV}$, the functions $\varphi(\lambda
t)/\varphi(t)$ and $\varphi(t)/\varphi(\lambda t)$ are bounded in
a neighborhood of $+\infty$. Therefore applying the theorem
\cite[Sec. 1.2]{Se85} on uniform convergence to a positive slowly
varying function $\chi_{1}$ such that
$\chi_{1}(\tau)\asymp\chi(\tau)$ with $\tau\gg1$, we can write
$$
\chi_{1}\left((\lambda t)^{\theta}\varphi(\lambda t)\right)\big/
\chi_{1}\left(t^{\theta}\varphi(t)\right)=
\chi_{1}\left(\frac{\lambda^{\theta}\varphi(\lambda
t)}{\varphi(t)}\:t^{\theta}\varphi(t)\right)\Big/
\chi_{1}\left(t^{\theta}\varphi(t)\right)\rightarrow1\quad\mbox{as}
\quad t\rightarrow+\infty.
$$
Here we use the limit $t^{\theta}\varphi(t)\rightarrow\infty$ as
$t\rightarrow+\infty$. Hence, the function
$\chi_{1}(t^{\theta}\varphi(t))$ is slowly varying at $+\infty$.
In addition,
$\chi(t^{\theta}\varphi(t))\asymp\chi_{1}(t^{\theta}\varphi(t))$
with $t\gg1$. Thus, the function $\chi(t^{\theta}\varphi(t))$
belongs to the set $\mathrm{QSV}$. Assertion (iv) is proved.
\end{proof}

\subsection{A refined scale over the Euclidean space.}
Let $n\in\mathbb{N}$. As usual, $\mathbb{R}^{n}$ denotes the
$n$-dimensional Euclidean space and $\mathcal{S}'(\mathbb{R}^{n})$
denotes the linear topological Schwartz space of tempered
distributions on $\mathbb{R}^{n}$. We use also the following
notations:
$\langle\xi\rangle=(1+\xi_{1}^{2}+\ldots+\xi_{n}^{2})^{1/2}$
denotes the smoothed modulus of a vector
$\xi=(\xi_{1},\ldots,\xi_{n})\in \mathbb{R}^{n}$ and $\widehat{u}$
denotes the Fourier transform of the distribution
$u\in\mathcal{S}'(\mathbb{R}^{n})$. We write an integral taken on
the space $\mathbb{R}^{n}$ without limits.

Let $\mathcal{M}$ denote the set of all functions
$\varphi:[1;+\infty)\rightarrow(0;+\infty)$ such that
\begin{description}
  \item[a] $\varphi$ is Borel measurable on the set $[1;+\infty)$;
  \item[b] functions $\varphi$ and $1/\varphi$ are bounded on
every closed interval $[1;b]$, where $1<b<+\infty$;
  \item[c] $\varphi\in\mathrm{QSV}$.
\end{description}

Let $s\in\mathbb{R}$, $\varphi\in\mathcal{M}$.

\begin{definition}
We denote by $H^{s,\varphi}(\mathbb{R}^{n})$ the space of all
distributions $u\in\mathcal{S}'(\mathbb{R}^{n})$ such that the
Fourier transform $\widehat{u}$ is a function locally Lebesgue
integrable on $\mathbb{R}^{n}$ which satisfies the inequality
$$
\int\langle\xi\rangle^{2s}\varphi^{2}(\langle\xi\rangle)\,|\widehat{u}(\xi)|^{2}\,
d\xi<\infty.
$$
The inner product in the space
$\mathrm{H}^{s,\varphi}(\mathbb{R}^{n})$ is defined by the formula
$$
(u,v)_{\mathrm{H}^{s,\varphi}(\mathbb{R}^{n})}:=
\int\langle\xi\rangle^{2s}\varphi^{2}(\langle\xi\rangle)
\,\widehat{u}(\xi)\,\overline{\widehat{v}(\xi)}\,d\xi
$$
and generates the norm in the usual way.
\end{definition}

The space $H^{s,\varphi}(\mathbb{R}^{n})$ is a special isotropic
Hilbert case of the spaces introduced by L.~H\"ormander
\cite[Sec.~2.2]{Her65}, \cite[Sec. 10.1]{Her86} and L. R.
Volevich, B. P. Paneah \cite[Sec. 2]{VoPa65}, \cite[Sec.
1.4.2]{Pa00}. Note that this space is actually defined with the
help of the function $\varphi_{s}(t)=t^{s}\varphi(t)$ regularly
varying  at $+\infty$ with the index $s$. However it is more
convenient for us to represent the parameter $\varphi_{s}$ as the
couple of two parameters $s$ and $\varphi$.

In the particular case where $\varphi\equiv1$ the space
$H^{s,\varphi}(\mathbb{R}^{n})$ coincides with the Sobolev space
$H^{s}(\mathbb{R}^{n})$. In general, the following inclusions are
true:
$$
\bigcup_{\varepsilon>0}H^{s+\varepsilon}(\mathbb{R}^{n})=:H^{s+}(\mathbb{R}^{n})
\subset H^{s,\varphi}(\mathbb{R}^{n})\subset
H^{s-}(\mathbb{R}^{n}):=\bigcap_{\varepsilon>0}H^{s-\varepsilon}(\mathbb{R}^{n}).
\leqno(3.2)
$$
They result from assertion (ii) of Theorem 3.3 and from the
definition of the set $\mathcal{M}$, according to which, for each
$\varepsilon>0$, there is a number $c_{\varepsilon}\geq1$ such
that $c_{\varepsilon}^{-1}t^{-\varepsilon}\leq\varphi(t)\leq
c_{\varepsilon}t^{\varepsilon}$ for $t\geq1$. Inclusions (3.2)
mean that, in the collection of the spaces
$$
\{H^{s,\varphi}(\mathbb{R}^{n}):s\in\mathbb{R},\varphi\in\mathcal{M}\,\},
\leqno(3.3)
$$
the function parameter $\varphi$ \textit{refines} the basic
(power) $s$-smoothness. Therefore it is natural to give the
following definition.

\vspace{0.3 cm}

\begin{definition}
The collection of function
spaces (3.3) is called a \textit{refined scale} over
$\mathbb{R}^{n}$ (with respect to the Sobolev scale).
\end{definition}

Besides the properties inherent to the H\"ormander spaces
\cite[Sec. 2.2 ]{Her65}, \cite[Sec. 10.1]{Her86} and the
Volevich-Paneah spaces \cite[Ch. I, II]{VoPa65}, \cite[Sec.
1.4]{Pa00}, the refined scale over $\mathbb{R}^{n}$ possesses the
following fundamental interpolation property:

\begin{theorem}
Let a function $\varphi\in\mathcal{M}$ and positive numbers
$\varepsilon,\delta$ be given. Let
$\psi(t):=t^{\,\varepsilon/(\varepsilon+\delta)}\,
\varphi(t^{1/(\varepsilon+\delta)})$ for $t\geq1$ and
$\psi(t):=\varphi(1)$ for $0<t<1$. Then the following assertions
are true:
\begin{itemize}
  \item[(i)] The function $\psi$ belongs to the set $\mathcal{B}$
and is an interpolation parameter.
  \item[(ii)] For an
arbitrary $s\in\mathbb{R}$, the equality of spaces
$$
\left[H^{s-\varepsilon}(\mathbb{R}^{n}),H^{s+\delta}(\mathbb{R}^{n})\right]_{\psi}
=H^{s,\varphi}(\mathbb{R}^{n})
$$
and equality of norms in them hold.
\end{itemize}
\end{theorem}

\begin{proof} \textit{Assertion} (i). By virtue of assertions
(ii), (iv) of Theorem 3.3, the function $\psi$ belongs to the set
$\mathcal{B}$ and is a function regular varying at $+\infty$ with
the index $\theta=\varepsilon/(\varepsilon+\delta)\in(0,\,1)$.
Therefore $\psi$ is an interpolation parameter because of Theorem
3.2. Assertion (i) is proved.

\textit{Assertion} (ii). Let $s\in\mathbb{R}$. It follows from the
properties of the Sobolev spaces that the couple
$[H^{s-\varepsilon}(\mathbb{R}^{n}),H^{s+\delta}(\mathbb{R}^{n})]$
is admissible and the pseudodifferential operator with symbol
$\langle\xi\rangle^{\varepsilon+\delta}$ is a generating operator
$J$ for this couple. Applying the Fourier transform
$\mathcal{F}:H^{s-\varepsilon}(\mathbb{R}^{n})\leftrightarrow
L_{2}(\mathbb{R}^{n},\langle\xi\rangle^{2(s-\varepsilon)}d\xi)$,
we reduce the operator $J$ to the form of multiplication by the
function $\langle\xi\rangle^{\varepsilon+\delta}$ of
$\xi\in\mathbb{R}^{n}$. Hence, the operator $\psi(J)$ is reduced
to the form of multiplication by the function
$\psi(\langle\xi\rangle^{\varepsilon+\delta})=
\langle\xi\rangle^{\varepsilon}\varphi(\langle\xi\rangle)$. This
permits us to write the following in view of (3.2):
$$
\left[H^{s-\varepsilon}(\mathbb{R}^{n}),H^{s+\delta}(\mathbb{R}^{n})\right]_{\psi}=
\left\{u\in H^{s-\varepsilon}(\mathbb{R}^{n}):
\langle\xi\rangle^{\varepsilon}\varphi(\langle\xi\rangle)\
\widehat{u}(\xi)\in
L_{2}(\mathbb{R}^{n},\langle\xi\rangle^{2(s-\varepsilon)}d\xi)\right\}
$$
$$
=\left\{u\in H^{s-\varepsilon}(\mathbb{R}^{n}):
\int\langle\xi\rangle^{2s}\varphi^{2}(\langle\xi\rangle)
\left|\widehat{u}(\xi)\right|^{2}d\xi<\infty\right\}=
H^{s-\varepsilon}(\mathbb{R}^{n})\cap
H^{s,\varphi}(\mathbb{R}^{n})=H^{s,\varphi}(\mathbb{R}^{n}).
$$
In addition the norm in the space
$\left[H^{s-\varepsilon}(\mathbb{R}^{n}),
H^{s+\delta}(\mathbb{R}^{n})\right]_{\psi}$ is equal to
$$
\|\psi(J)u\|_{H^{s-\varepsilon}(\mathbb{R}^{n})}=
\left(\int|\langle\xi\rangle^{\varepsilon}\varphi(\langle\xi\rangle)\
\widehat{u}(\xi)|^{2}\,\langle\xi\rangle^{2(s-\varepsilon)}\
d\xi\right)^{1/2}=\|u\|_{H^{s,\varphi}(\mathbb{R}^{n})}.
$$
Assertion (ii) is proved.
\end{proof}

\subsection{A refined scale over a closed manifold}
Further let $\Gamma$ be a closed (that is compact and without a
boundary) infinitely smooth manifold of dimension $n$. We suppose
that a certain $C^{\infty}$-density $dx$ is defined on $\Gamma$.
We denote by $\mathcal{D}'(\Gamma)$ the linear topological space
of all distributions on $\Gamma$. Thus $\mathcal{D}'(\Gamma)$ is
the space antidual to the space $C^{\infty}(\Gamma)$ with respect
to the natural extension of the scalar product in
$L_{2}(\Gamma,dx)=:L_{2}(\Gamma)$ by continuity. This extension is
denoted by $(f,w)_{\Gamma}$ for $f\in\mathcal{D}'(\Gamma)$, $w\in
C^{\infty}(\Gamma)$.

The refined scale over the manifold $\Gamma$ is constructed from
scale (3.3) in the following way. We choose a finite atlas from
the $C^{\infty}$-structure on $\Gamma$ consisting of the local
charts $\alpha_{j}:\mathbb{R}^{n}\leftrightarrow U_{j}$,
$j=1,\ldots,r$. Here the open sets $U_{j}$ form the finite
covering of the manifold $\Gamma$. Let functions $\chi_{j}\in
C^{\infty}(\Gamma)$, $j=1,\ldots,r$, form a partition of unity on
$\Gamma$ satisfying the condition $\mathrm{supp}\,\chi_{j}\subset
U_{j}$. As before, $s\in\mathbb{R}, \varphi\in\mathcal{M}$.

\begin{definition}
We denote by $H^{s,\varphi}(\Gamma)$ the space of all
distributions $f\in\mathcal{D}'(\Gamma)$ such that
$(\chi_{j}f)\circ\alpha_{j}\in H^{s,\varphi}(\mathbb{R}^{n})$ for
each $j=1,\ldots,r$. Here $(\chi_{j}f)\circ\alpha_{j}$ is the
representation of the distribution $\chi_{j}f$ in the local chart
$\alpha_{j}$. The inner product in the space
$H^{s,\varphi}(\Gamma)$ is defined by the formula
$$
(f,g)_{H^{s,\varphi}(\Gamma)}:=\sum_{j=1}^{r}\,((\chi_{j}f)\circ\alpha_{j},
(\chi_{j}\,g)\circ\alpha_{j})_{H^{s,\varphi}(\mathbb{R}^{n})}
$$
and induces the norm in the usual way.
\end{definition}

\begin{definition}
The collection of function
spaces
$\{H^{s,\varphi}(\Gamma):s\in\mathbb{R},\varphi\in\mathcal{M}\}$
is called a \textit{refined scale} over the closed manifold
$\Gamma$.
\end{definition}

In the particular case where $\varphi\equiv1$ the space
$H^{s,\varphi}(\Gamma)$ coincides with the Sobolev space
$H^{s}(\Gamma)$. Sobolev spaces are known \cite[Sec. 2.6]{Her65},
\cite[Sec. 7.5]{Shubin78} to be complete and independent (up to
equivalence of norms) of the choice of the atlas and the partition
of unity. We will show that every space $H^{s,\varphi}(\Gamma)$
can be obtained by the interpolation of the proper couple of
Sobolev's spaces. It implies that the space
$H^{s,\varphi}(\Gamma)$ is Hilbert and independent of this choice.

\begin{theorem}
Let a function
$\varphi\in\mathcal{M}$ and positive numbers $\varepsilon,\delta$
be given. Then, for each $s\in\mathbb{R}$, the equality of spaces
$$
\left[H^{s-\varepsilon}(\Gamma),\
H^{s+\delta}(\Gamma)\right]_{\psi}=H^{s,\varphi}(\Gamma)\quad\mbox{with
equivalence of norms} \leqno(3.4)
$$
hold. Here $\psi$ is the interpolation parameter from Theorem
$3.4$.
\end{theorem}

\begin{proof}
The couple of the Sobolev spaces on the left-hand side of equality
(3.4) is admissible \cite[Sec. 7.5, 7.6]{Shubin78}. We deduce this
equality from Theorem 3.4 with the help of the well known method
of "rectification"\, and "sewing"\, of the manifold $\Gamma$.
According to Definition 3.5, the linear mapping of "rectification"
$$
T:\,f\mapsto(\,(\chi_{1}f)\circ\alpha_{1},\ldots,(\chi_{r}f)\circ\alpha_{r}\,),
\quad f\in\mathcal{D}'(\Gamma),
$$
defines the isometric operators
$$
T:\,H^{\sigma}(\Gamma)\rightarrow(H^{\sigma}(\mathbb{R}^{n}))^{r},\quad
\sigma\in\mathbb{R}, \leqno(3.5)
$$
$$
T:\,H^{s,\varphi}(\Gamma)\rightarrow(H^{s,\varphi}(\mathbb{R}^{n}))^{r}.
\leqno(3.6)
$$
Since $\psi$ is the interpolation parameter and operators (3.5)
are bounded for $\sigma\in\{s-\varepsilon,s+\delta\}$, the bounded
operator
$$
T:\,\left[H^{s-\varepsilon}(\Gamma),H^{s+\delta}(\Gamma)\right]_{\psi}
\rightarrow\left[\,(H^{s-\varepsilon}(\mathbb{R}^{n}))^{r},\,
(H^{s+\delta}(\mathbb{R}^{n}))^{r}\,\right]_{\psi} \leqno(3.7)
$$ exists.
By virtue of Theorems 2.5, 3.4, the following equalities of spaces
and norms in them are true:
$$
\left[\,(H^{s-\varepsilon}(\mathbb{R}^{n}))^{r},\,
(H^{s+\delta}(\mathbb{R}^{n}))^{r}\,\right]_{\psi}=
\left(\,\left[H^{s-\varepsilon}(\mathbb{R}^{n}),
H^{s+\delta}(\mathbb{R}^{n})\right]_{\psi}\,\right)^{r}
=(H^{s,\varphi}(\mathbb{R}^{n}))^{r}. \leqno(3.8)
$$
Thus, since operator (3.7) is bounded, so is the operator
$$
T:\,\left[H^{s-\varepsilon}(\Gamma),H^{s+\delta}(\Gamma)\right]_{\psi}
\rightarrow(H^{s,\varphi}(\mathbb{R}^{n}))^{r}. \leqno(3.9)
$$

Now we construct for $T$ the left inverse operator $K$ of
"sewing"\, of the manifold $\Gamma$. For each $j=1,\ldots,r$ we
take a function $\eta_{j}\in C_{0}^{\infty}(\mathbb{R}^{n})$ such
that $\eta_{j}=1$ on the set
$\alpha_{j}^{-1}(\mathrm{supp}\,\chi_{j})$. Let us consider the
linear mapping
$$
K:\,(h_{1},\ldots,h_{r})\mapsto\sum_{j=1}^{r}\,
\Theta_{j}\left((\eta_{j}h_{j})\circ\alpha_{j}^{-1}\right),\quad
h_{1},\ldots,h_{r}\in\mathcal{S}'(\mathbb{R}^{n}).
$$
Here $(\eta_{j}h_{j})\circ\alpha_{j}^{-1}$ is a distribution in
the open set $U_{j}\subseteq\Gamma$ such that its representative
in the local map $\alpha_{j}$ has the form $\eta_{j}h_{j}$. In
addition, $\Theta_{j}$ denotes the operator of extension by zero
from $U_{j}$ to $\Gamma$. This operator is well defined on
distributions with support belonging to $U_{j}$. By the choice of
the functions $\chi_{j}$, $\eta_{j}$, we have
$$
KTf=\sum_{j=1}^{r}\,\Theta_{j}\left((\eta_{j}\,
((\chi_{j}f)\circ\alpha_{j}))\circ\alpha_{j}^{-1}\right)=
\sum_{j=1}^{r}\,\Theta_{j}\left(\,
(\chi_{j}f)\circ\alpha_{j}\circ\alpha_{j}^{-1}\right)=
\sum_{j=1}^{r}\,\chi_{j}f=f,
$$
that is
$$
KTf=f\quad\mbox{for each}\quad f\in\mathcal{D}'(\Gamma).
\leqno(3.10)
$$

Let us show that the linear mapping $K$ defines the bounded
operator
$$
K:\,(H^{s,\varphi}(\mathbb{R}^{n}))^{r}\rightarrow
H^{s,\varphi}(\Gamma).\leqno(3.11)
$$
For an arbitrary vector $h=(h_{1},\ldots,h_{r})$ from the space
$(H^{s,\varphi}(\mathbb{R}^{n}))^{r}$, we write
$$
\bigl\|Kh\bigr\|^{2}_{H^{s,\varphi}(\Gamma)}
=\sum_{l=1}^{r}\;\bigl\|(\chi_{l}\,Kh)
\circ\alpha_{\,l}\bigr\|_{H^{s,\varphi}(\mathbb{R}^{n})}^{2}
=\sum_{l=1}^{r}\,\Bigl\|\Bigl(\chi_{\,l}\,\sum_{j=1}^{r}\
\Theta_{j}\bigl((\eta_{j}h_{j})\circ\alpha_{j}^{-1}\bigr)\Bigr)
\circ\alpha_{\,l}\,\Bigr\|_{H^{\,s,\varphi}(\mathbb{R}^{n})}^{2}
$$
$$
=\sum_{l=1}^{r}\;\Bigl\|\,\sum_{j=1}^{r}(\eta_{j,l}\,h_{j})
\circ\beta_{\,j,l}\,\Bigr\|_{H^{s,\varphi}(\mathbb{R}^{n})}^{2}
\leq\sum_{l=1}^{r}\;\Bigl(\sum_{j=1}^{r}\,\bigl\|(\eta_{j,l}\,h_{j})\circ
\beta_{j,l}\,\bigr\|_{H^{\,s,\varphi}(\mathbb{R}^{n})}\Bigr)^{2}.\leqno(3.12)
$$
Here $\eta_{j,l}:=(\chi_{\,l}\circ\alpha_{j})\,\eta_{j}\in
C_{0}^{\infty}(\mathbb{R}^{n})$ and
$\beta_{j,l}:\,\mathbb{R}^{n}\,\leftrightarrow\,\mathbb{R}^{n}$ is
a $C^{\infty}$-diffeomorphism such that
$\beta_{j,l}\,=\alpha_{j}^{-1}\circ\alpha_{l}$ in a neighborhood
of $\mathrm{supp}\,\eta_{j,l}$ and $\beta_{j,l}(x)=x$ for all
$x\in\mathbb{R}^{n}$ sufficiently large in modulus. The operator
of multiplication by a function of the class
$C_{0}^{\infty}(\mathbb{R}^{n})$ and the operator of change of
variables $u\mapsto u\circ\beta_{j,l}$ are known \cite[Theorems
B.1.7, B.1.8]{Her87} to be  bounded on every space
$H^{\sigma}(\mathbb{R}^{n})$ with $\sigma\in\mathbb{R}$. Therefore
the linear operator $v\mapsto (\eta_{j,l}\,v)\circ\beta_{j,l}$ is
bounded on the space $H^{\sigma}(\mathbb{R}^{n})$. Then
boundedness of this operator in the space
$H^{s,\varphi}(\mathbb{R}^{n})$ follows by Theorem 3.3. Hence
relations (3.12) imply the estimate
$$
\bigl\|Kh\bigr\|_{H^{s,\varphi}(\Gamma)}^{2}\leq
c\,\sum_{j=1}^{r}\
\bigl\|h_{j}\bigr\|_{H^{s,\varphi}(\mathbb{R}^{n})}^{2},
$$
where the constant $c>0$ is independent of
$h=(h_{1},\ldots,h_{r})$. Thus, operator (3.11) is bounded for
each $s\in\mathbb{R}$, $\varphi\in\mathcal{M}$.

In particular, the operators
$K:\,(H^{\sigma}(\mathbb{R}^{n}))^{r}\rightarrow
H^{\sigma}(\Gamma)$ with $\sigma\in\mathbb{R}$ are bounded. Let us
set $\sigma\in\{s-\varepsilon,\,s+\delta\}$ and use the
interpolation with the parameter $\psi$. Due to equality (3.8), we
obtain the bounded operator
$$
K:\,(H^{s,\varphi}(\mathbb{R}^{n}))^{r}\rightarrow\
\left[H^{s-\varepsilon}(\Gamma),\,H^{s+\delta}(\Gamma)\right]_{\psi}.\leqno(3.13)
$$
Now formulae (3.6), (3.13) and (3.10) imply that the identity
mapping $KT$ establishes the continuous imbedding of the space
$H^{s,\varphi}(\Gamma)$ into the interpolation space
$[H^{s-\varepsilon}(\Gamma),\,H^{s+\delta}(\Gamma)]_{\psi}$.
Moreover, formulae (3.10) and (3.13) imply that the same mapping
$KT$ establishes also the inverse continuous imbedding. Theorem
3.5 is proved.
\end{proof}

The following properties of the refined scale over the manifold
$\Gamma$ can be deduced from Theorem 3.5 and the interpolation
properties established in Section 2.

\begin{theorem}
Let $s\in\mathbb{R}$ and $\varphi,\varphi_{1}\in\mathcal{M}$. The
following assertions are true:
\begin{itemize}
  \item [(i)] The space $H^{s,\varphi}(\Gamma)$ is
Hilbert separable and does not depend (up to equivalence of norms)
on the choice of an atlas for $\Gamma$ and partition of unity used
in Definition $3.5$.
  \item [(ii)] The set $C^{\infty}(\Gamma)$ is dense in the space
$H^{s,\varphi}(\Gamma)$.
  \item [(iii)] For each
$\varepsilon>0$, the compact and dense imbedding
$H^{s+\varepsilon,\varphi_{1}}(\Gamma)\hookrightarrow
H^{s,\varphi}(\Gamma)$ holds.
  \item [(iv)] Suppose that the function $\varphi/\varphi_{1}$
is bounded in a neighborhood of $+\infty$. Then continuous dense
imbedding $H^{s+\varepsilon,\varphi_{1}}(\Gamma)\hookrightarrow
H^{s,\varphi}(\Gamma)$ is valid. It is compact if
$\varphi(t)/\varphi_{1}(t)\rightarrow0$ as $t\rightarrow+\infty$.
  \item [(v)] The spaces $H^{s,\varphi}(\Gamma)$ and
$H^{-s,1/\varphi}(\Gamma)$ are mutually dual (up to equivalence of
norms) with respect to the extension of the inner product in
$L_{2}(\Gamma)$ by continuity.
  \end{itemize}
\end{theorem}

\begin{proof}
\textit{Assertion} (i). The space $H^{s,\varphi}(\Gamma)$ is
Hilbert and separable because, according to Theorem 3.5, this
space is obtained by the interpolation of a certain couple of the
Sobolev spaces. Let us consider two couples $\mathcal{A}_{1}$ and
$\mathcal{A}_{2}$ each of which consists of an atlas and a
partition of unity on $\Gamma$. We denote by
$H^{s,\varphi}(\Gamma,\mathcal{A}_{j})$ and
$H^{\sigma}(\Gamma,\mathcal{A}_{j})$ respectively the spaces from
the refined scale and the Sobolev spaces which correspond to the
couple $\mathcal{A}_{j}$, where $j=1,\,2$. For the Sobolev spaces,
the identity mapping establishes the topological isomorphism
$I:H^{\sigma}(\Gamma,\mathcal{A}_{1})\leftrightarrow
H^{\sigma}(\Gamma,\mathcal{A}_{2})$ for each
$\sigma\in\mathbb{R}$. Let us set $\sigma=s\mp1$ and use the
interpolation with the parameter $\psi$ from Theorem 3.4. By
Theorem 3.5 we arrive at the topological isomorphism
$I:H^{s,\varphi}(\Gamma,\mathcal{A}_{1})\leftrightarrow
H^{s,\varphi}(\Gamma,\mathcal{A}_{2})$. It means that the space
$H^{s,\varphi}(\Gamma)$ is independent of the choice of an atlas
and a unity partition mentioned above. Assertion (i) is proved.

\textit{Assertion} (ii). By virtue of Theorems 2.1 and 3.5, we
have the continuous dense imbedding
$H^{s+\delta}(\Gamma)\hookrightarrow H^{s,\varphi}(\Gamma)$.
Besides, the set $C^{\infty}(\Gamma)$ is dense in the Sobolev
space  $H^{s+\delta}(\Gamma)$ \cite[Proposition 7.4]{Shubin78}.
These two facts imply assertion (ii).

\textit{Assertion} (iii). Assume that $\varepsilon>0$. By Theorem
3.5, there  exist interpolation parameters
$\chi,\eta\in\mathcal{B}$ such that the following equalities of
spaces with equivalence of norms in them is true:
$$
\left[H^{s+\varepsilon/2}(\Gamma),H^{s+2\varepsilon}(\Gamma)\right]_{\chi}=
H^{s+\varepsilon,\varphi_{1}}(\Gamma)\quad\mbox{and}\quad
\left[H^{s-\varepsilon}(\Gamma),H^{s+\varepsilon/3}(\Gamma)\right]_{\eta}=
H^{s,\varphi}(\Gamma).
$$
It implies by Theorem 2.1 the next chain of continuous imbeddings
$$
H^{s+\varepsilon,\varphi_{1}}(\Gamma)\hookrightarrow
H^{s+\varepsilon/2}(\Gamma)\hookrightarrow
H^{s+\varepsilon/3}(\Gamma)\hookrightarrow H^{s,\varphi}(\Gamma).
$$
Here the central imbedding of Sobolev spaces is compact
\cite[Theorem 7.4]{Shubin78}. Therefore the imbedding
$H^{s+\varepsilon,\varphi_{1}}(\Gamma)\hookrightarrow
H^{s,\varphi}(\Gamma)$ is compact as well. This imbedding is dense
because of assertion (ii). Assertion (iii) is proved.

\textit{Assertion} (iv). Let us assume that the function
$\varphi/\varphi_{1}$ is bounded in a neighborhood of $+\infty$.
By Theorem 3.5, we have the following equalities of spaces with
equivalence of norms in them:
$$
\left[H^{s-1}(\Gamma),H^{s+1}(\Gamma)\right]_{\psi}
=H^{s,\varphi}(\Gamma)\quad\mbox{and}\quad
\left[H^{s-1}(\Gamma),H^{s+1}(\Gamma)\right]_{\psi_{1}}
=H^{s,\varphi_{1}}(\Gamma).
$$
Here the interpolation parameters $\psi,\psi_{1}\in\mathcal{B}$
satisfy the condition
$$
\psi(t)/\psi_{1}(t)=\varphi(t^{1/2})/\varphi_{1}(t^{1/2})\quad
\mbox{for}\quad t\geq1.
$$
Hence, the function $\psi/\psi_{1}$ is bounded in a neighborhood
of $+\infty$ that, by Theorem 2.2, implies the continuous dense
imbedding $H^{s,\varphi_{1}}(\Gamma)\hookrightarrow
H^{s,\varphi}(\Gamma)$. Now suppose that
$\varphi(t)/\varphi_{1}(t)\rightarrow0$ as $t\rightarrow+\infty$.
It implies the limit $\psi(t)/\psi_{1}(t)\rightarrow0$ as
$t\rightarrow+\infty$. In addition, we recall that the imbedding
of Sobolev spaces $H^{s+1}(\Gamma)\hookrightarrow H^{s-1}(\Gamma)$
is compact. It follows by Theorem 2.2 that the imbedding
$H^{s,\varphi_{1}}(\Gamma)\hookrightarrow H^{s,\varphi}(\Gamma)$
is compact as well. Assertion (iv) is proved.

\textit{Assertion} (v) is known (see e.g. \cite[Theorem
7.7]{Shubin78}) in the case $\varphi\equiv1$. From this the case
of an arbitrary $\varphi\in\mathcal{M}$ can be obtained as
follows. First we note that $1/\varphi\in\mathcal{M}$ and
therefore the space $H^{-s,1/\varphi}(\Gamma)$ is well defined.
The Sobolev spaces $H^{s\pm1}(\Gamma)$ and $H^{-s\mp1}(\Gamma)$
are mutually dual with respect to the extension of the inner
product in $L_{2}(\Gamma)$ by continuity. This means that the
linear mapping $Q:w\mapsto(\,\cdot\,,\overline{w})_{\Gamma}$,
$w\in C^{\infty}(\Gamma)$, is extended by continuity to the
topological isomorphisms
$Q:H^{s\mp1}(\Gamma)\leftrightarrow(H^{-s\pm1}(\Gamma))'$. Let us
apply to them the interpolation with the parameter $\psi$ from
Theorem 3.5 in the case where $\varepsilon=\delta=1$. We obtain
one more topological isomorphism
$$
Q:\left[H^{s-1}(\Gamma),H^{s+1}(\Gamma)\right]_{\psi}\leftrightarrow
\left[(H^{-s+1}(\Gamma))',(H^{-s-1}(\Gamma))'\right]_{\psi}.
\leqno(3.14)
$$
Here the left-hand interpolation space equals to
$H^{s,\varphi}(\Gamma)$ and, by Theorem 2.4, the right-hand one
can be written as
$$
\left[(H^{-s+1}(\Gamma))',(H^{-s-1}(\Gamma))'\right]_{\psi}=
\left[H^{-s-1}(\Gamma),H^{-s+1}(\Gamma)\right]_{\chi}'=
(H^{-s,1/\varphi}(\Gamma))'.
$$
Let us note that the last equality is valid because
$\chi(t):=t/\psi(t)=t^{1/2}/\varphi(t^{1/2})$ for $t\geq1$. Thus,
(3.14) implies the topological isomorphism
$Q:H^{s,\varphi}(\Gamma)\leftrightarrow(H^{-s,1/\varphi}(\Gamma))'$,
which means the mutual duality of the spaces
$H^{s,\varphi}(\Gamma)$ and $H^{-s,1/\varphi}(\Gamma)$ in the
sense mentioned above. Assertion (v) is proved.
\end{proof}

The refined scale is closed with respect to the interpolation with
a function parameter regular varying at $+\infty$.

\begin{theorem}
Let $s_{0},s_{1}\in\mathbb{R}$, $s_{0}\leq s_{1}$ and
$\varphi_{0},\varphi_{1}\in\mathcal{M}$. In the case where
$s_{0}=s_{1}$ we suppose that the function
$\varphi_{0}/\varphi_{1}$ is bounded in a neighborhood of
$+\infty$. Let $\psi\in\mathcal{B}$ be a function regularly
varying at $+\infty$ with the index $\theta$, where $0<\theta<1$.
By Theorem $3.2$, $\psi$ is an interpolation parameter. We
represent it as $\psi(t)=t^{\theta}\chi(t)$ with
$\chi\in\mathrm{QSV}$. Let us set $s:=(1-\theta)s_{0}+\theta
s_{1}$ and
$$
\varphi(t):=\varphi_{0}^{1-\theta}(t)\,\varphi_{1}^{\theta}(t)\,
\chi\left(t^{s_{1}-s_{0}}\varphi_{1}(t)/\varphi_{0}(t)\right)\quad\mbox{for}\quad
t\geq1.
$$
Then $\varphi\in\mathcal{M}$ and
$$
\left[\,H^{s_{0},\varphi_{0}}(\Gamma),
H^{s_{1},\varphi_{1}}(\Gamma)\,\right]_{\psi}=
H^{s,\varphi}(\Gamma)\quad\mbox{with equivalence of norms}.
$$
\end{theorem}

\begin{proof}
This theorem is a direct consequence of
Theorems 3.5 and 2.3.
\end{proof}

\begin{remark}
Theorem 3.7 is true in the limiting case where
$\theta=0$ or $\theta=1$ under additional supposition that the
function $\psi$ is quasiconcave in a neighborhood of $+\infty$.
Then, by Theorem\;2.7, $\psi$ is an interpolation parameter. For
example, Theorem 3.7 is true for each of the functions
$\psi(t):=\ln^{r}t$ and $\psi(t):=t/\ln^{r} t$, where $t\gg1$ and
$r>0$.
\end{remark}

\subsection{An alternative definition of the refined scale.} Let
$A$ be an elliptic pseudodifferential operator on $\Gamma$ of
order $m>0$. We suppose that the operator
$A:C^{\infty}(\Gamma)\rightarrow C^{\infty}(\Gamma)$ is positive
on the space $L_{2}(\Gamma)$, that is there exists a number $r>0$
such that
$$
(Lu,u)_{\Gamma}\geq r\,(u,u)_{\Gamma}\quad\mbox{for each}\quad
u\in C^{\infty}(\Gamma). \leqno(3.15)
$$
In the present subsection, $(\,\cdot\,,\,\cdot\,)_{\Gamma}$ is the
inner product in $L_{2}(\Gamma)$.

We denote by $A_{0}$ the closure of the operator
$A:C^{\infty}(\Gamma)\rightarrow C^{\infty}(\Gamma)$ on the space
$L_{2}(\Gamma)$. This closure exists and has the domain
$H^{m}(\Gamma)$ because the operator $A$ is elliptic on $\Gamma$
\cite[Corollary 8.3]{Shubin78}, \cite[Theorem 2.3.5]{Agr90}. The
pseudodifferential operator $A$ is formally self-adjoint due to
condition (3.15). Hence \cite[Theorem 8.3]{Shubin78},
\cite[Theorem 2.3.7]{Agr90}, $A_{0}$ is an unbounded self-adjoint
operator in the space $L_{2}(\Gamma)$ with
$\mathrm{Spec}\,A_{0}\subseteq[r,+\infty)$. In particular, we have
$0\notin\mathrm{Spec}\,A_{0}$, that implies the topological
isomorphism
$$
A:\,H^{s+m,\varphi}(\Gamma)\,\leftrightarrow\,H^{s,\varphi}(\Gamma)
\quad\mbox{for each}\quad s\in\mathbb{R},\;\varphi\in\mathcal{M}.
\leqno(3.16)
$$
In the Sobolev case where $\varphi\equiv1$ this result is well
known (see e.g. \cite[Theorem 8.1, Proposition 8.5]{Shubin78},
\cite[Theorem 19.2.1]{Her87}, \cite[Sec. 2.3]{Agr90}). The general
case of an arbitrary $\varphi\in\mathcal{M}$ follows immediately
from the case $\varphi\equiv1$ by virtue of Theorem 3.5.

Let $s\in\mathbb{R}$ and $\varphi\in\mathcal{M}$. We set
$\varphi_{s}(t):=t^{s/m}\varphi(t^{1/m})$ for $t\geq1$ and,
moreover, $\varphi_{s}(t):=\varphi(1)$ for $0<t<1$. Since the
function $\varphi$ is positive and Borel measurable on the
semiaxis $(0,+\infty)$, a self-adjoint operator
$\varphi_{s}(A_{0})$ is well-defined in the space $L_{2}(\Gamma)$
as the function $\varphi$ of $A_{0}$.

\begin{lemma}
The following assertions are true.
\begin{itemize}
  \item [(i)] The domain of the operator $\varphi_{s}(A_{0})$
contains the set $C^{\infty}(\Gamma)$.
  \item [(ii)] The mapping
$$
f\,\mapsto\,\|\varphi_{s}(A_{0})f\|_{L_{2}(\Gamma)},\quad f\in
C^{\infty}(\Gamma), \leqno(3.17)
$$
is a norm in the space $C^{\infty}(\Gamma)$.
  \end{itemize}
\end{lemma}

\begin{proof} \textit{Assertion} (i). Let us choose an integer $k$ so that
$k>s/m$. Since $\varphi\in\mathcal{M}$, the function $\varphi_{s}$
is bounded on every compact subset of the semiaxis $(0,+\infty)$
and, moreover, $t^{-k}\varphi_{s}(t)\rightarrow0$ as
$t\rightarrow+\infty$ because of assertions (ii), (iv) of Theorem
3.3. Hence, there is a number $c>0$ such that $\varphi_{s}(t)\leq
c\,t^{k}$ for $t\geq r$. Let us consider the unbounded operator
$A_{0}^{k}$ on the space $L_{2}(\Gamma)$. Since
$A:C^{\infty}(\Gamma)\rightarrow C^{\infty}(\Gamma)$, we can write
$C^{\infty}(\Gamma)\subset\mathrm{Dom}\,A_{0}^{k}\subset
\mathrm{Dom}\,\varphi_{s}(A_{0})$. Assertion (i) is proved.

\textit{Assertion} (ii). According to assertion (i), mapping
(3.17) is well-defined. For this mapping, all norm properties  are
evident except for the positive definiteness property. Let us
prove it. Applying the spectral theorem, we can write for an
arbitrary function $f\in C^{\infty}(\Gamma)$:
$$
\bigl\|\varphi_{s}(A_{0})f\bigr\|_{L_{2}(\Gamma)}^{2}=
\int_{r}^{+\infty}\varphi_{s}^{2}(t)\,d(E_{t}f,f)_{\Gamma}\quad\mbox{and}
\quad\bigl\|f\bigr\|_{L_{2}(\Gamma)}^{2}=
\int_{r}^{+\infty}d(E_{t}f,f)_{\Gamma}. \leqno(3.18)
$$
Here $E_{t}$, $t\geq r$, is the resolution of identity in the
space $L_{2}(\Gamma)$ which corresponds to the self-adjoint
operator $A_{0}$. If
$\|\varphi_{s}(A_{0})f\|_{L_{2}(\Gamma)}^{2}=0$, then from the
first equality in (3.18) and from the inequality $\varphi_{s}>0$
it follows that the measure $(E(\cdot)f,f)_{\Gamma}$ of the set
$[r,+\infty)$ is equal to 0. Now the second equality in (3.18)
implies that $f=0$ on $\Gamma$. Assertion (ii) is proved.
\end{proof}

\begin{definition}
The space
$H^{s,\varphi}_{A}(\Gamma)$ is a completion of the space
$C^{\infty}(\Gamma)$ with respect to norm (3.17).
\end{definition}

The space $H^{s,\varphi}_{A}(\Gamma)$ is Hilbert one because norm
(3.17) is generated by the inner product
$(\varphi_{s}(A_{0})f,\varphi_{s}(A_{0})g)_{\Gamma}$ of functions
$f,g\in C^{\infty}(\Gamma)$.

\vspace{0.3 cm}

\begin{theorem}
For arbitrary $s\in\mathbb{R}$, $\varphi\in\mathcal{M}$, the norms
in the spaces $H^{s,\varphi}_{A}(\Gamma)$ and
$H^{s,\varphi}(\Gamma)$ are equivalent on the dense linear
manifold $C^{\infty}(\Gamma)$. Thus,
$H^{s,\varphi}_{A}(\Gamma)=H^{s,\varphi}(\Gamma)$ up to
equivalence of norms.
\end{theorem}

\begin{proof}
At first let us suppose that $s>0$. We choose $k\in\mathbb{N}$ so
that $k\,m>s$. Since the operator $A_{0}^{k}$ is closed and
positive on $L_{2}(\Gamma)$, its domain $\mathrm{Dom}\,A_{0}^{k}$
is Hilbert space with respect to the inner product
$(A_{0}^{k}f,A_{0}^{k}g)_{\Gamma}$ of functions $f,g$. Note that
the couple of spaces $[L_{2}(\Gamma),\mathrm{Dom}\,A_{0}^{k}]$ is
admissible, and the operator $A_{0}^{k}$ is a generating one for
it. Moreover, since $A_{0}^{k}$ is a closure of the elliptic
pseudodifferential operator $A^{k}$ on $L_{2}(\Gamma)$, the spaces
$\mathrm{Dom}\,A_{0}^{k}$ and $H^{km}(\Gamma)$ are equal up to
equivalent norms. Let a function $\psi$ be the interpolation
parameter from Theorems 3.3, 3.4 with $\varepsilon=s$ and
$\delta=k\,m-s$. Then $\psi(t^{k})=\varphi_{s}(t)$ for $t>0$, so
by Theorem 3.4 we can write
$$
\bigl\|f\bigr\|_{H^{s,\varphi}(\Gamma)}\asymp
\bigl\|f\bigr\|_{[H^{0}(\Gamma),\,H^{km}(\Gamma)]_{\psi}}\asymp
\bigl\|f\bigr\|_{[L_{2}(\Gamma),\,\mathrm{Dom}\,A_{0}^{k\,}]_{\psi}}=
\bigl\|\psi(A_{0}^{k})f\bigr\|_{L_{2}(\Gamma)}=
\bigl\|\varphi_{s}(A_{0})f\bigr\|_{L_{2}(\Gamma)},
$$
for each $f\in C^{\infty}(\Gamma)$.

Now let the number $s\in\mathbb{R}$ be arbitrary. Choose
$k\in\mathbb{N}$ so that $s+k\,m>0$. As has been proved,
$$
\bigl\|g\bigr\|_{H^{s+km,\varphi}(\Gamma)}\asymp
\bigl\|\varphi_{s+km}(A_{0})\,g\bigr\|_{L_{2}(\Gamma)},\quad g\in
C^{\infty}(\Gamma). \leqno(3.19)
$$
The following topological isomorphism holds due to (3.16) :
$$
A^{k}:\,H^{\sigma+km,\varphi}(\Gamma)\leftrightarrow
H^{\sigma,\varphi}(\Gamma)\quad\mbox{for
each}\quad\sigma\in\mathbb{R}. \leqno(3.20)
$$
Denote by $A^{-k}$ the inverse operator to $A^{k}$. For every
function $f\in C^{\infty}(\Gamma)$, we have $A^{-k}f\in
C^{\infty}(\Gamma)$ and $A_{0}^{k}A^{-k}f=f$. Hence, by virtue of
(3.20), (3.19), we can write
$$
\bigl\|f\bigr\|_{H^{s,\varphi}(\Gamma)}\asymp
\bigl\|A^{-k}f\bigr\|_{H^{s+km,\varphi}(\Gamma)}\asymp
\bigl\|\varphi_{s+km}(A_{0})A^{-k}f\bigr\|_{L_{2}(\Gamma)}
$$
$$
=\bigl\|\varphi_{s}(A_{0})A_{0}^{k}A^{-k}f\bigr\|_{L_{2}(\Gamma)}=
\bigl\|\varphi_{s}(A_{0})f\bigr\|_{L_{2}(\Gamma)},\quad f\in
C^{\infty}(\Gamma).
$$
Theorem 3.8 is proved.
\end{proof}

\begin{theorem}
Let $s\geq0$ and $\varphi\in\mathcal{M}$. In the case where $s=0$
we suppose that the function $1/\varphi$ is bounded in a
neighborhood of $+\infty$. Then the space $H^{s,\varphi}(\Gamma)$
coincides with the domain of the operator $\varphi_{s}(A_{0})$ and
the norm in the space $H^{s,\varphi}(\Gamma)$ is equivalent to the
graph norm of the operator $\varphi_{s}(A_{0})$.
\end{theorem}

\begin{proof}
The domain $\mathrm{Dom}\,\varphi_{s}(A_{0})$ of the closed
operator $\varphi_{s}(A_{0})$ is Hilbert space with respect to the
scalar product of the graph of this operator. Let us prove that
the norms in the spaces $\mathrm{Dom}\,\varphi_{s}(A_{0})$ and
$H^{s,\varphi}_{A}(\Gamma)$ are equivalent on the dense linear
manifold $C^{\infty}(\Gamma)$. By Theorem 3.8, it will imply the
present theorem. According to the condition of the present theorem
and by virtue of assertion (ii) of Theorem 3.3, there is a number
$c>0$ such that $\varphi_{s}(t)\geq c$ for $t>0$. Therefore
$$
\bigl\|\varphi_{s}(A_{0})f\bigr\|_{L_{2}(\Gamma)}\geq c\,
\bigl\|f\bigr\|_{L_{2}(\Gamma)}\quad\mbox{for each}\quad f\in
C^{\infty}(\Gamma).
$$
It yields the equivalence of norms mentioned above. It remains to
prove the density of the set $C^{\infty}(\Gamma)$ in the space
$\mathrm{Dom}\,\varphi_{s}(A_{0})$.

Let $f\in\mathrm{Dom}\,\varphi_{s}(A_{0})$. Since
$\varphi_{s}(A_{0})f\in L_{2}(\Gamma)$, there is a sequence of
functions $h_{j}\in C^{\infty}(\Gamma)$ such that
$h_{j}\rightarrow\varphi_{s}(A_{0})f$ in $L_{2}(\Gamma)$ as
$j\rightarrow\infty$. Note that the operator
$\varphi_{s}^{-1}(A_{0})$ is bounded on the space $L_{2}(\Gamma)$
because $1/\varphi_{s}(t)\leq 1/c$ for $t>0$. Hence,
$$
f_{j}:=\varphi_{s}^{-1}(A_{0})h_{j}\rightarrow
f\quad\mbox{and}\quad
\varphi_{s}(A_{0})f_{j}=h_{j}\rightarrow\varphi_{s}(A_{0})f\quad\mbox{in}\quad
L_{2}(\Gamma)\quad\mbox{as}\quad j\rightarrow\infty.
$$
In other words, $f_{j}\rightarrow f$ with respect to the graph
norm of the operator $\varphi_{s}(A_{0})$. Moreover, since
$h_{j}\in C^{\infty}(\Gamma)$, then
$f_{j}=A_{0}^{-k}\varphi_{s}^{-1}(A_{0})A_{0}^{k}\,h_{j}\in
H^{km}(\Gamma)$ for every $k\in\mathbb{N}$. Consequently,
$f_{j}\in C^{\infty}(\Gamma)$ and the density of the set
$C^{\infty}(\Gamma)$ in the space
$\mathrm{Dom}\,\varphi_{s}(A_{0})$ is established. Theorem 3.9 is
proved.
\end{proof}

A significant example of the operator $A$ investigated above is
the operator $1-\triangle_{\Gamma}$, where $\triangle_{\Gamma}$ is
the Beltrami-Laplace operator on the Riemannian manifold $\Gamma$
(then $m=2$).

\end{document}